\documentclass[preprint, floatsintext]{elsarticle}
\usepackage[english]{babel}
\usepackage[utf8]{inputenc}
\usepackage[displaymath, mathlines]{lineno}
\usepackage{amsfonts}
\usepackage{amsmath}
\usepackage{amsthm}
\usepackage{bm}
\usepackage{booktabs}
\usepackage[labelsep=period]{caption}
\usepackage{enumitem}
\usepackage{fancyhdr}
\usepackage{float}
\usepackage{geometry}
\usepackage{graphicx}
\usepackage{lastpage}
\usepackage{listings}
\usepackage{longtable}
\usepackage{subcaption}
\usepackage{units}
\usepackage[svgnames]{xcolor}

\newtheorem{theorem}{Theorem}[section]

\newtheorem{lemma}[theorem]{Lemma}
\newtheorem{definition}[theorem]{Definition}

\newcommand{\realR}{\mathbb{R}}

\newcommand{\bhat}[1]{\hat{\bm{#1}}}
\newcommand{\dint}[1]{~\text{d}#1}

\newcommand{\domain}{\Omega}
\newcommand{\CartDomain}{\hat{\Omega}}
\newcommand{\embeddedRegion}{\Omega_E}

\newcommand{\refElement}{\hat{D}}
\newcommand{\curve}{\mathcal{C}}

\newcommand{\RHSc}{A}
\newcommand{\SRDc}{S}
\newcommand{\RHSd}{\bm{\mathrm{A}}}

\newcommand{\uSysC}{\bm{U}}

\newcommand{\bu}{\bm{u}}

\newcommand{\bV}{\bm{V}}

\newcommand{\grad}{\nabla}

\renewcommand{\div}{\grad \cdot}

\newcommand{\pd}[3]{\frac{\partial^{#3} #1}{\partial#2^{#3}}}

\newcommand{\Grad} {\ensuremath{\nabla}}
\newcommand{\Div} {\ensuremath{\nabla\cdot}}

\newcommand{\jump}[1] {\ensuremath{\LRs{\![#1]\!}}}

\newcommand{\LRp}[1]{\left( #1 \right)}
\newcommand{\LRs}[1]{\left[ #1 \right]}

\newcommand{\eval}[2][\right]{\relax
  \ifx#1\right\relax \left.\fi#2#1\rvert}

\begin{document} % =============================================================

\begin{frontmatter}

% Author Information -----------------------------------------------------------

\cortext[cor1]{Corresponding author(s)}

\author[1]{Christina G. Taylor\corref{cor1}}
\ead{cgt3@rice.edu}

\author[2]{Lucas C. Wilcox}

\author[1]{Jesse Chan\corref{cor1}}
\ead{jesse.chan@rice.edu}

\affiliation[1]{organization={Department of Computational Applied Mathematics and Operations Research, Rice University}}
	
\affiliation[2]{organization={Department of Applied Mathematics,
	Naval Postgraduate School}}

% Paper Information ------------------------------------------------------------
% Old: Energy Stable State Redistribution Cut Cell Discont. Galerkin Methods for Wave Propagation
\title{An Energy Stable High-Order Cut Cell Discontinuous Galerkin Method with State Redistribution for Wave Propagation}
\date{Spring 2024}


\begin{abstract}

Cut meshes are a type of mesh that is formed by allowing embedded boundaries to ``cut" a simple underlying mesh resulting in a hybrid mesh of cut and standard elements. While cut meshes can allow complex boundaries to be represented well regardless of the mesh resolution, their arbitrarily shaped and sized cut elements can present issues such as the \textit{small cell problem}, where small cut elements can result in a severely restricted CFL condition. State redistribution, a technique developed by Berger and Giuliani \cite{berger-stateRedistr}, can be used to address the small cell problem. In this work, we pair state redistribution with a high-order discontinuous Galerkin scheme that is $L_2$ energy stable for arbitrary quadrature. We prove that state redistribution can be added to a provably $L_2$ energy stable discontinuous Galerkin method on a cut mesh without damaging the scheme's $L_2$ stability. We numerically verify the high order accuracy and stability of our scheme on two-dimensional wave propagation problems. 
\end{abstract}

\end{frontmatter}

\section{Introduction} %==============================================================================
\begin{figure}[H]
\centering
\begin{subfigure}{0.49\textwidth}
\centering
\includegraphics[width=\textwidth]{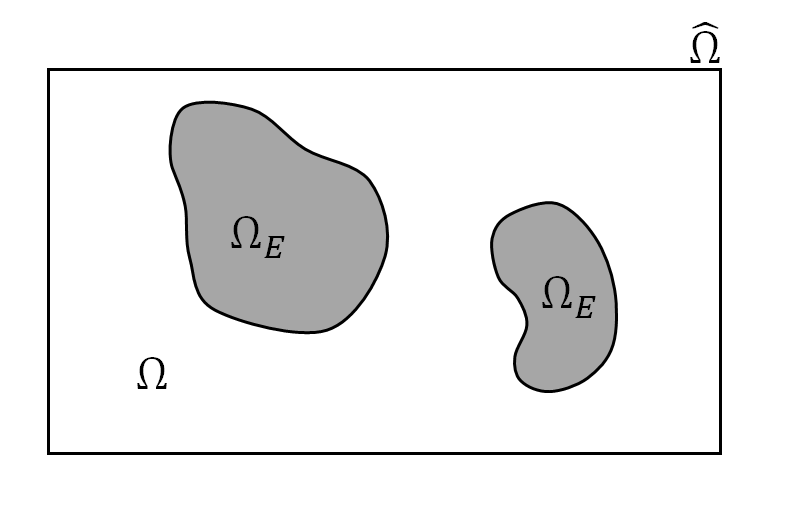}
\caption{Example domain. }\label{fig:domainRef}
\end{subfigure}
\begin{subfigure}{0.49\textwidth}
\centering
\includegraphics[width=\textwidth]{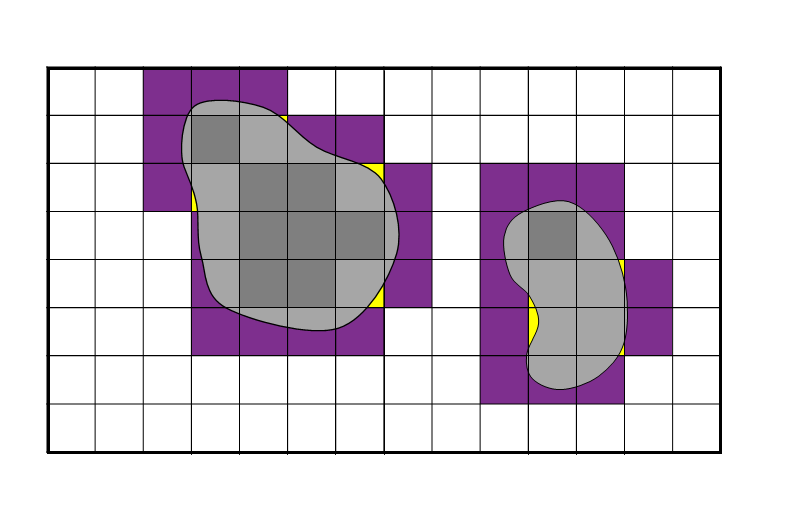}
\caption{Example cut mesh.}\label{fig:cutMeshRef}
\end{subfigure}

\caption{An example of an embedding Cartesian domain, $\hat{\Omega}$, with embedded boundaries defining regions $\Omega_{E}$ to exclude and the resulting simulation domain $\Omega$. Once meshed, elements remain Cartesian (white), become cut cells (active portion purple and yellow), or get excluded from the simulation domain (dark grey). Note cut meshes can result in extremely small or skewed elements, such as those shown in yellow.}\label{fig:domainMesh-ref}
\end{figure}

Here we present a provably energy stable discontinuous Galerkin (DG) method on 2D cut meshes using state redistribution. Our method falls underneath the umbrella of Embedded Boundary methods, where embedded regions $\embeddedRegion$ are removed from the overall embedding domain, $\CartDomain$, which we take to be Cartesian, to yield the final PDE domain 
\begin{equation}
\domain = \CartDomain ~\backslash~ \embeddedRegion.
\end{equation} 
Here we limit ourselves only to 2D domains and thus $\CartDomain, \embeddedRegion, \domain \in \realR^2$. An example of a domain with embedded boundaries is shown in Figure \ref{fig:domainRef}.

Of interest to us is the simulation of hyperbolic conservation laws. Such conservation laws govern a wide variety of natural phenomena and include the wave equation, shallow water equations, and compressible Euler equations. The inclusion of embedded boundaries allows simulation of fluid flow around objects and more geometrically complex domains. These laws feature wave-like behavior, making the acoustic wave equation, 
\begin{align}
\frac{1}{c^2}\pd{p}{t}{} &+ \Div \bm{u} = 0, \\[12pt]
\pd{\bm{u}}{t}{} &+ \Grad p = 0,
\end{align}
a model problem before considering more complex hyperbolic conservation laws \cite{dafermos-HyperbolicConservationLawsBook}. Here we limit our focus to the scalar wave equation, but our scheme can be applied to any linear symmetric hyperbolic system \cite{warburton-skewSymm}. 

In general, linear hyperbolic conservation laws satisfy a continuous $L_2$ energy conservation property: the $L_2$ energy is constant in time. When dealing with linear conservation laws, certain numerical methods preserve this energy conservation at the semi-discrete level and are thus said to be energy conservative. If the numerical method additionally is able to dissipate energy, i.e. the $L_2$ energy can decrease in time, it is said to be energy stable. For high-order methods, energy stable formulations are particularly pertinent as high-order methods can suffer from instability \cite{wang-highVsLowOrder-HOinstability}. High-order methods are of interest to us as they have greater accuracy per degrees of freedom \cite{wang-highVsLowOrder-HOinstability} and are better suited to unsteady flows \cite{visbal-highOrder-unsteadyFlows} due to their lower dissipation compared to low-order methods \cite{ainsworth-highOrder-dispersionError}. However, high-order methods also feature more degrees freedom per element and thus are more expensive under mesh refinement. To avoid the need to refine our mesh we appeal to cut meshes to capture embedded boundaries. 

Cut meshes are a form of hybrid mesh and can be linear or curvilinear in nature. In cut meshes the (typically Cartesian) embedding domain is discretized using standard meshing techniques without regard for the embedded boundaries to create a background mesh. Then, the embedded boundaries are allowed to ``cut" and exclude elements from the background mesh resulting in a hybrid mesh of cut and standard elements. An example is given in Figure \ref{fig:cutMeshRef}.

Cut meshes are attractive for use with embedded boundaries as they are relatively easy to generate and can capture complex geometries well while maintaining computational efficiency \cite{berger-cutStabilityAccuracy}. Cut meshes have a long history in computational fluid dynamics field with one of their first appearances in \cite{reed-originalCutCell} in the 1970s. Since their introduction they have been used for a number of applications such as fluid dynamics \cite{quirk-EarlyCutCell, ingram-cutCellDevelopments, chung-cellMerging}, fluid-structure interaction \cite{schott-cutFSIMonolithic}, front tracking \cite{udaykumar-phaseFront}, moving boundaries \cite{udaykumar-movingCutCell}, and physiology models \cite{berre-rognes}. Here we pair our high-order DG scheme with a Cartesian cut mesh that features a coarse, computationally efficient Cartesian quadrilateral background mesh with curvilinear cut elements. 

Despite their advantages, cut meshes, which features arbitrarily shaped/sized cut elements, can pose difficulties for many methods and classes of problems. A recurring challenge is the construction of volume quadrature rules on cut elements, though this challenge is relative straightforward to address using techniques from the curvilinear meshing community \cite{warburton-skewSymm, chan-skewSymm} when exact quadrature is not needed. A greater challenge when simulating hyperbolic problems on cut meshes is the \textit{small cell problem}, where the arbitrarily shaped (and sized) cut elements severely restrict the CFL (Courant-Friedichs-Lewy; i.e., the maximum stable time step) condition of the simulation.

\subsection{Methods for Addressing the Small Cell Problem}
Addressing the small cell problem is an active area of research for which a number of methods have been developed. One way to avoid the small cell problem altogether is to use implicit time integration as in \cite{xie-ImplicitCut}. However, the cost and complexity of implicit time integration can be prohibitive. Hybrid explicit-implicit time integration schemes such as in  \cite{may-explicitImplicit-original, may-explicitImplicit-accuracy} seek to marry the stability of implicit time integration with the practicality of explicit time integration by using implicit integration on small elements and explicit integration elsewhere. Such local time integration schemes can also be fully explicit by using subcycling \cite{schoeder-localTimestepping} on small cells. However the use of localized time integration requires coupling of different time integration schemes and may not sufficiently stabilize all small elements \cite{schoeder-localTimestepping}.

More standard method-of-lines approaches using special DG and finite element (FE) schemes have also been developed that make use of special penalty terms or ghost values to stabilize small elements. One such example is presented in \cite{sticko-stabilizedBilinearForm} where penalization of jumps in normal derivatives is used to build stabilization into the spatial discretization's mass and stiffness matrices. This has been done for DG and finite element \cite{fu-cutDGPenaltyTerms, gurkan-aPriori-geomIndep, sticko-nistche-aPrioriError, sticko-stabilizedBilinearForm} formulations. The use of penalty terms can yield the attractive result of being able to assert energy stability \cite{sticko-stabilizedBilinearForm} or provide a priori error estimates such as in \cite{sticko-nistche-aPrioriError, gurkan-aPriori-geomIndep}. However, these error estimates may not always be robust for all geometries \cite{dePrenter-noteOnNistche} though some provably are \cite{gurkan-aPriori-geomIndep}.

Another penalty-term based method is the domain-of-dependence method, which was developed solely for hyperbolic problems and uses penalty terms to redistribute mass from a cut cell to its neighbors \cite{engwer-DoD-original}. While it was originally developed on triangular meshes for the advection equation \cite{engwer-DoD-original} it has since been extended to general linear hyperbolic problems on general elements \cite{birke-DoD-linearHyperbolic} and in 1D to nonlinear hyperbolic problems \cite{may-DoD-nonlinearHyperbolic}. It has been shown to be $L_2$ energy stable \cite{birke-DoD-linearHyperbolic}, total variation diminishing in the $L_1$ norm \cite{engwer-DoD-original}, and monotone for piecewise constant solutions \cite{streitburger-DoD-monotonicity}. 

While penalty terms seek to provide stabilization by modifying the DG or FE scheme's formulation, flux-based stabilization is also possible for finite volume (FV) and DG methods. One such flux-based stabilization approach is $h$-box methods, which seek to solve a Riemann problem over a ``box" of length/scale $h$ at cell interfaces to compute the numerical fluxes for each cut cell. $h$-box methods were first developed for 1D FV cut cells in \cite{berger-hBox-original}. This method was further developed in \cite{helzel-hBox, berger-hBox-cut} for adaptation to 2D and to improve accuracy and computational efficiency \cite{berger-simplifiedhBox}. Another flux-based stabilization method is flux redistribution, which was originally developed for front tracking in a hybrid finite difference, finite volume scheme in \cite{chern-fluxRedistrbution-early}. Similar to the $h$-box method, flux redistribution seeks to adjust small cut elements/cells' flux values but typically does so using flux values from neighboring elements. While flux redistribution is conservative, it is only first-order accurate \cite{chern-fluxRedistrbution-early, colella-fluxRedistribution}.

In the same spirit as flux redistribution is cell merging, also called cell agglomeration, where small elements are merged with neighbors to eliminate small elements altogether. Cell merging is a common solution to the small cell problem \cite{chung-cellMerging, berger-cutStabilityAccuracy} and can be used in conjunction with other stabilization methods \cite{schoeder-localTimestepping}. However, determining a provably stable time step with cell merging can be difficult in 2D and 3D meshes \cite{berger-CellMerging}. Similar to cell merging is cell linking, which has the benefit of maintaining the original, unmerged mesh \cite{muralidharan-cellLinking, kirkpatrick-cellLinking, cecere-cellLinking}.

To address the small cell problem in our solver we use \textit{state redistribution}, a method developed by Berger and Giuliani in \cite{berger-stateRedistr}. State redistribution is in a sense an extension of cell merging and linking but with the issue of order dependency resolved. Importantly, it is high-order accurate and conservative in the sense that is preserves the average of the solution and polynomial degree \cite{berger-stateRedistr}. 

In this paper we prove that when state redistribution is applied to an energy stable DG scheme the resulting scheme remains energy stable. The proof of energy stability under state redistribution is this project's main contribution to the Embedded Boundary method literature. 

The rest of this paper is organized as follows: in Section \ref{sec:methods} we describe our DG formulation, which is energy stable for arbitrary quadrature rules on cut meshes, before describing state redistribution in more detail in preparation for Section \ref{sec:proof}, where we present our proof of energy stability under state redistribution. Next, in Section \ref{sec:results} we numerically verify the theoretical properties of our scheme and compare our solver to the benchmark ``Pacman" model for 2D acoustic wave propagation \cite{zeigelwanger-Pacman}. Lastly in the Appendices we highlight some of the implementational details of our scheme, such as our approach to representing the embedded boundaries and constructing volume quadrature on cut elements.

\section{Methods}\label{sec:methods} %===================================================================================

\subsection{The DG Formulation}
For our DG solver we consider the linear acoustic wave equation. On a given domain $\Omega$ with appropriate initial and boundary conditions the acoustic wave equation is given by
\begin{align}
\frac{1}{c^2}\pd{p}{t}{} &+ \Div \bm{u} = 0,\label{eqn:dpdt} \\[12pt]
\pd{\bm{u}}{t}{} &+ \Grad p = 0,\label{eqn:dudt}
\end{align}
where $c$ is the speed of sound in the medium, $p$ is the pressure, and $\bu \in \realR^2$ is the velocity. We assume a Cartesian embedding domain $\CartDomain$ over which we impose a Cartesian rectangular mesh of $n_x \times n_y$ elements. As stated previously our true domain $\domain$ is related to the Cartesian domain $\CartDomain$ by 
\begin{equation}
\domain = \CartDomain ~\backslash~  \embeddedRegion,
\end{equation}
where $\embeddedRegion$ is the set of embedded regions to be excluded as illustrated previously in Figure \ref{fig:domainMesh-ref}. The domain $\domain$ is then decomposed into a cut mesh of $N_h$ non-overlapping Cartesian and cut elements, $D^k$, satisfying
\begin{equation}
    \domain = \bigcup_{k=1}^{N_h} D^k.
\end{equation}

For a given polynomial degree $N$, we define the solution on Cartesian elements using a reference element $\refElement = [-1,1]^2 $. For the reference element solution space we use $\mathcal{Q}^N(\hat{D})$, the space of tensor products of degree $N$ polynomials in the $\hat{x}$ and $\hat{y}$ direction on the reference element. For this space we use a tensor product basis of 1D Lagrange polynomials,
\begin{equation}
    \phi_{ij}(\hat{x},\hat{y}) = \ell_i(\hat{x})\ell_j(\hat{y}), \quad i,j = 0,1,2,..., N.
\end{equation}

We define the solution on cut elements similarly with a nodal basis, but due to the arbitrary shape of cut elements we can no longer exploit a tensor product structure or use a reference element. On cut elements we use the space of total degree $N$ polynomials $\mathbb{P}^{N}(D^k) = \text{span}\{ x^i y^j ~|~ i + j \le N, i, j \in \mathbb{Z}^+_0 \}$. For convenience, we define these polynomials nodally on the physical cut element using 2D Lagrange polynomials,
\begin{equation}
    \varphi_{i}(x,y) = \ell_{i}(x,y), \quad \ell_i(x_j, y_j) = \delta_{ij} \text{~for~} 0 \le i,j,\quad i + j \le N,
\end{equation}
whose nodal points are approximate Fekete points \cite{sommariva-feketeQR}. Fekete points are a class of interpolation nodes that yield well-conditioned interpolation matrices \cite{sommariva-feketeQR} with the added benefit of being easy to generate; however, any set of unisolvent points will suffice. Our cut quadrature rules also use approximate Fekete points, but require a larger set of nodes which changes the points' distribution. In Appendix \ref{sec:VolumeQuadr} we elaborate on how we generate approximate Fekete points for use as volume quadrature nodes via a moment-fitting scheme. The same process is used to generate the polynomial basis nodes and is a form of moment fitting \cite{davis-momentFitting}.

We define our solution and test space as $\mathcal{Q}^N(D^k)$ on Cartesian elements and $\mathbb{P}^N(D^k)$ on cut elements. Given our solution $(p, \bm{u})$ and test functions $(q, \bm{v})$, the standard DG formulation \cite{hesthaven-TheDGBook} for the acoustic wave equation on a single element $D^k$ is given by
 \begin{align}
\int_{D^k}\frac{1}{c^2}\pd{p}{t}{} q = \int_{D^k} \Div \bm{u} q + \frac{1}{2}\int_{\partial D^k} \jump{\bm{u}}\cdot\bm{n} q,\label{eqn:stdFormulation-p}\\[10pt]
\int_{D^k}\pd{\bm{u}}{t}{} \cdot\bm{v} = \int_{D^k} \Grad p \cdot \bm{v}  + \frac{1}{2}\int_{\partial D^k} \jump{p}\bm{n}\cdot\bm{v}.\label{eqn:stdFormulation-u}
\end{align}
where $\bm{n}$ is the outward normal of the element $D^k$ and $\jump{\cdot}$ is used to denote the ``jump" between the solution on $D^k$ and its neighbors
\begin{equation}
    \jump{p} = p^+ - p, \quad\quad \jump{\bm{u}} = \bm{u}^+ - \bm{u}.
\end{equation}
where the ``+" superscript is used to denote the exterior value of a quantity on a neighboring element.

Under exact quadrature, this formulation is energy stable on affine meshes and standard DG approximation spaces.  However, cut elements necessitate custom quadrature rules, which in our current implementation do not provide exactness guarantees. To maintain energy stability we convert the standard formulation to its skew-symmetric form which is energy stable under arbitrary quadrature \cite{warburton-skewSymm, chan-skewSymm}. To derive the skew symmetric formulation, we split the volume terms into halves 
\begin{align}
    \int_{D^k} \Div \bm{u} q &= \frac{1}{2}\int_{D^k} \Div \bm{u} q + \frac{1}{2}\int_{D^k} \Div \bm{u} q,\\
    \int_{D^k} \Grad p \cdot \bm{v} &= \frac{1}{2}\int_{D^k} \Grad p \cdot \bm{v} + \frac{1}{2}\int_{D^k} \Grad p \cdot \bm{v},
\end{align}
and then apply integration by parts to \textit{one} of the halved terms. This allows terms to be combined/canceled with the surface integral terms to yield a new expression for the RHS (right-hand side):
\begin{align}
    \int_{D^k} \Div \bm{u} q + \frac{1}{2}\int_{\partial D^k} \jump{\bm{u}} \cdot{\bm{n}} q = \frac{1}{2}\int_{D^k} (\Div \bm{u} q - \bm{u}\cdot \Grad q) + \frac{1}{2}\int_{\partial D^k} \bm{u}^+\cdot{\bm{n}}q,\\[10pt]
    \int_{D^k} \Grad p \cdot \bm{v}  + \frac{1}{2}\int_{\partial D^k} \jump{p}\bm{n}\bm{v} = \frac{1}{2}\int_{D^k} (\Grad p \cdot \bm{v} - p\Div \bm{v})  + \frac{1}{2}\int_{\partial D^k} p^+ \bm{n} \cdot \bm{v}.
\end{align}
Substituting these expansions into the original formulation in Equations \eqref{eqn:stdFormulation-p} and \eqref{eqn:stdFormulation-u} yields the skew-symmetric formulation
\begin{align}
    \int_{D^k}\frac{1}{c^2}\pd{p}{t}{} q = \frac{1}{2}\int_{D^k} (\Div \bm{u} q - \bm{u}\cdot \Grad q) + \frac{1}{2}\int_{\partial D^k} \bm{u}^+\cdot{\bm{n}}q,\\[14pt]
    \int_{D^k}\pd{\bm{u}}{t}{}\cdot \bm{v} = \frac{1}{2}\int_{D^k} (\Grad p \cdot \bm{v} - p\Div \bm{v})  + \frac{1}{2}\int_{\partial D^k} p^+ \bm{n}\cdot \bm{v}.
\end{align}
Lastly we add upwinding penalty terms to the skew-symmetric formulation. These terms dissipate energy via penalizing jumps in the solution at element boundaries. Adding these terms yields our final formulation:
\begin{align}
    \int_{D^k}\frac{1}{c^2}\pd{p}{t}{} q = \frac{1}{2}\int_{D^k} (\Div \bm{u} q - \bm{u}\cdot \Grad q) + \frac{1}{2}\int_{\partial D^k} (\bm{u}^+\cdot{\bm{n}}q + \tau_p\jump{p}q)\label{eqn:wave-finalFormulation-p},\\[14pt]
    \int_{D^k}\pd{\bm{u}}{t}{} \cdot\bm{v} = \frac{1}{2}\int_{D^k} (\Grad p \cdot \bm{v} - p\Div \bm{v})  + \frac{1}{2}\int_{\partial D^k} (p^+\bm{n}\cdot\bm{v} + \tau_u\jump{\bm{u}}\cdot \bm{v}).  \label{eqn:wave-finalFormulation-u}
\end{align}
In all of our simulations we set the penalty parameters $\tau_p$ and $\tau_u$ equal to \nicefrac{1}{2}; in general they must satisfy $\tau_p, \tau_u \ge 0$. By adding energy dissipation, the penalty terms renders our final formulation energy stable \cite{warburton-skewSymm, chan-ESCurvilinear}. While we assume here that the speed of sound in the medium $c$ is constant for simplicity's sake, varying $c$ can also be accommodated by incorporating its effect into the mass matrix as shown in \cite{chan-heterogenMedia}.

\subsection{State Redistribution}

Next we briefly describe the action of state redistribution. State redistribution stabilizes the solution on small elements via a careful merging and redistribution of the solution with the solution on neighboring elements. In subsequent work since its original publication, a newer, upgraded version of state redistribution, called weighted state redistribution, uses a different weighting for the projections to allow state redistribution's effect to be smoothly activated \cite{giuliani-weightedSRD}. Here we use the original, volume-weighted version of state redistribution for high-order DG methods.

In state redistribution, we first identify cut elements in need of stabilization. A common condition is cut elements with volume less than one-half the volume of the background Cartesian element, which we also use. For each small element $D^k$, we determine a \textit{merging neighborhood}, over which we merge the constituent elements' solutions via a volume-weighted projection, $\Pi_k$. Mathematically we denote the merge neighborhood as 
\begin{equation}
M_k = \{k, k_1, k_2, ...\},
\end{equation}
the set of all element indices in element $D^k$'s merge neighborhood. For elements that do not require merging $M_k = \{k\}.$

\begin{figure}[H]
\centering
   \begin{subfigure}[b]{0.75\textwidth}
        \centering
        \includegraphics[width=\textwidth]{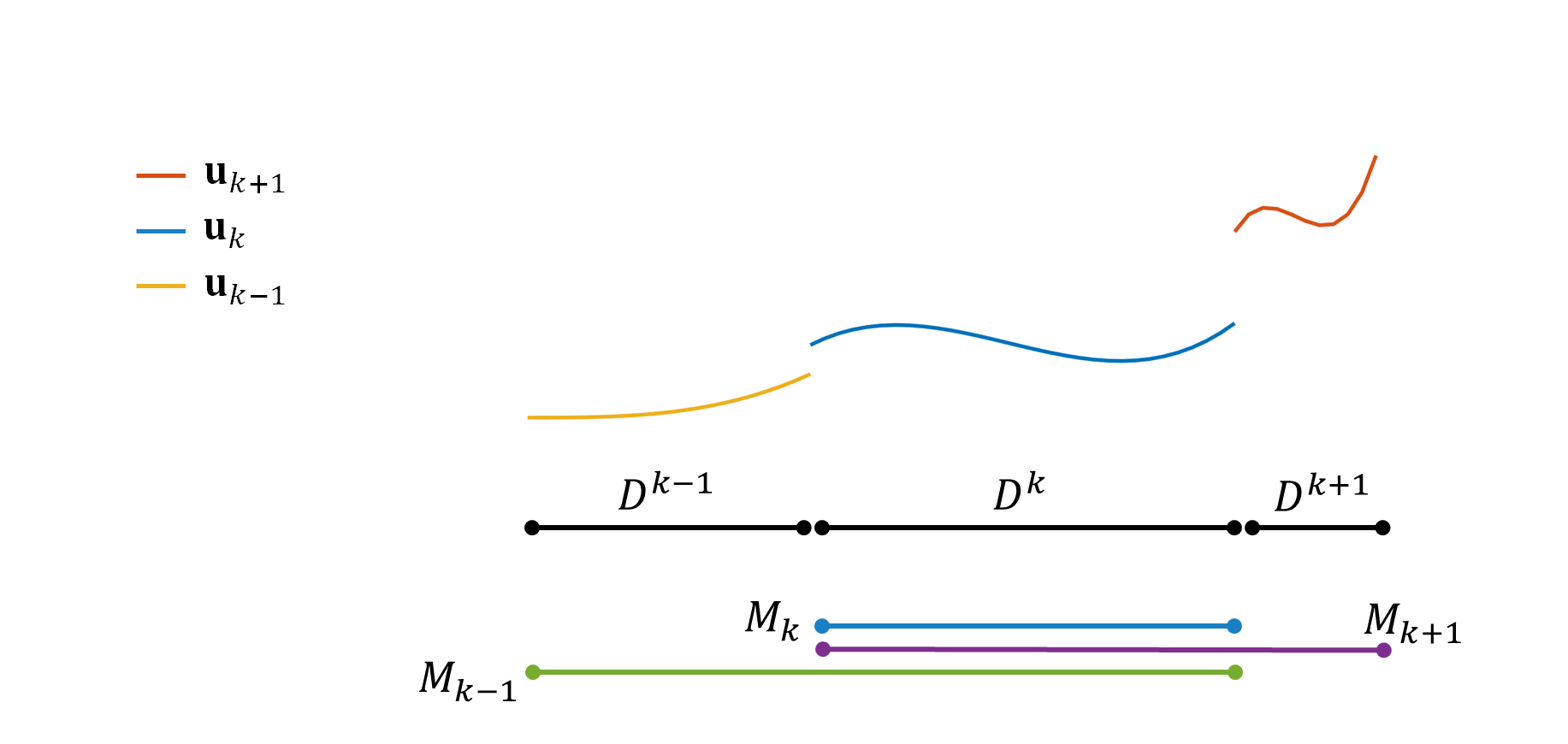}
        \caption{Here we denote the solution on element $D^k$ as $\bu_k$. Elements $D^{(i-1)}$ and $D^{(i+1)}$ are in need of merging. For each element, including full-sized element $D^k$, we compute merge neighborhoods $M_{k-1}, M_{k}, M_{k-1}$.}
   \end{subfigure}

   \begin{subfigure}[b]{0.75\textwidth}
        \centering
        \includegraphics[width=\textwidth]{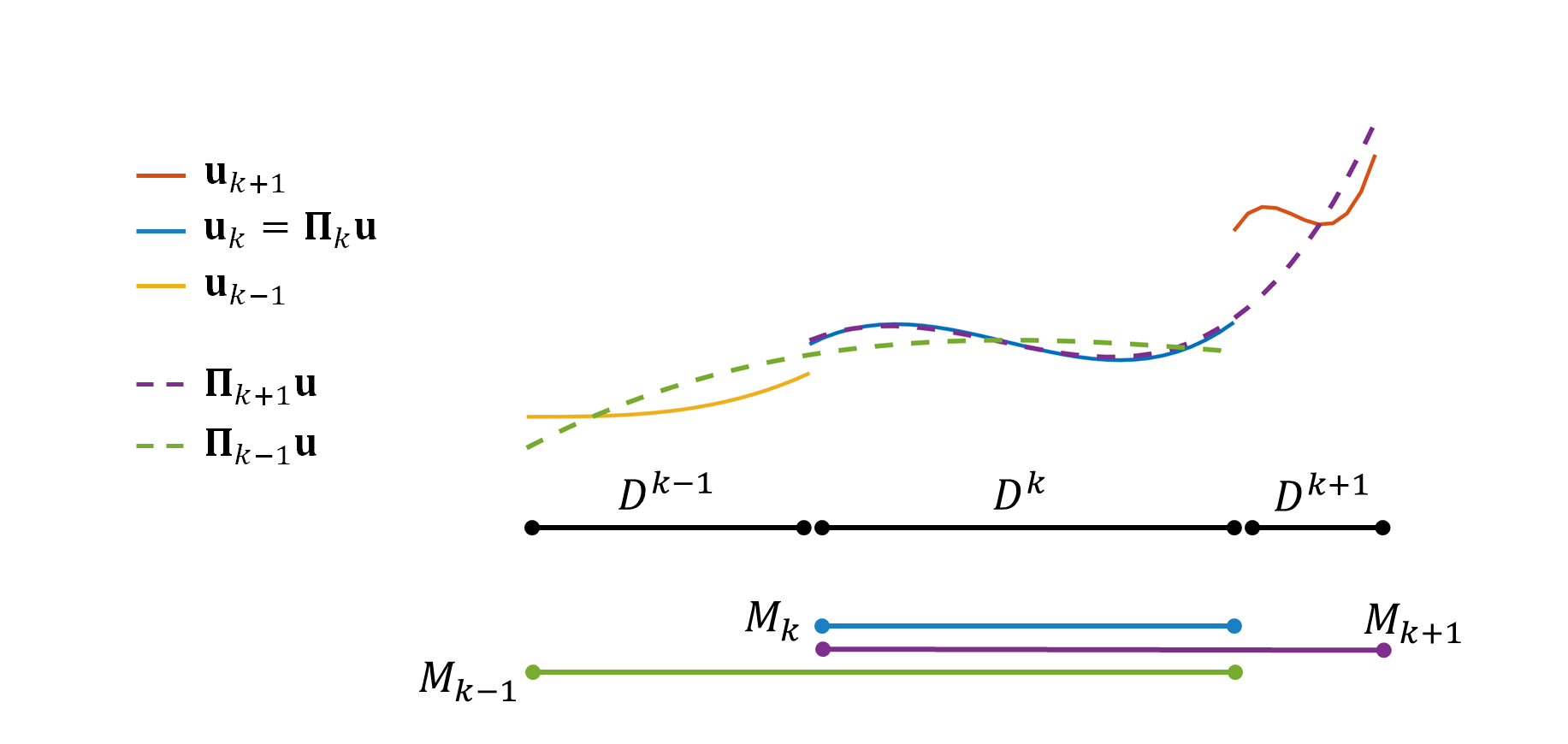}
        \caption{Over each merge neighborhood we compute the projected solutions. Notice in the case of a full-sized element the projected solution on its merge neighborhood is merely the original solution on that element.}
   \end{subfigure}

   \begin{subfigure}[b]{0.75\textwidth}
        \centering
        \includegraphics[width=\textwidth]{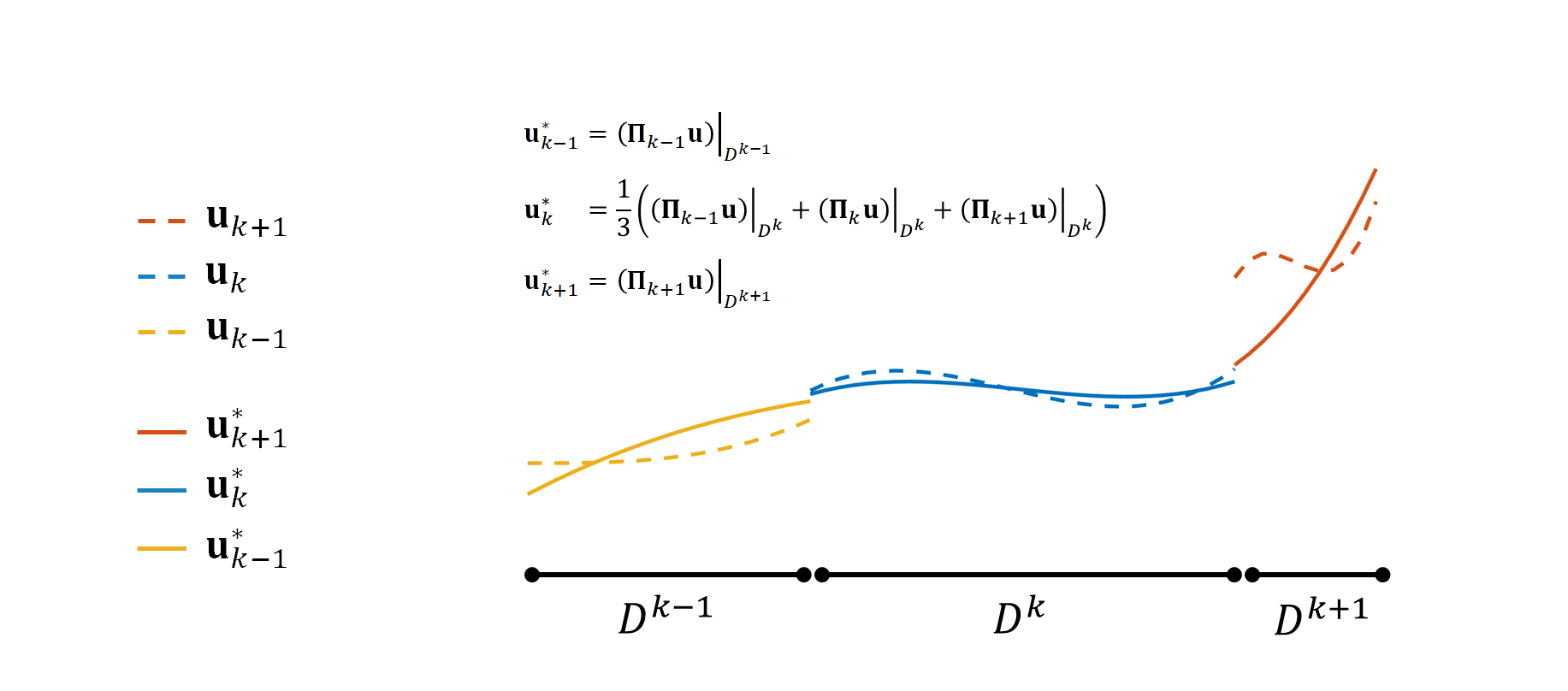}
        \caption{To compute the final updated solution $\bm{u}^*$ on each element we take the average of all projected solutions (restricted to the element of interest) to which that element contributed.}
   \end{subfigure}

\caption{The three-step process of state redistribution in 1D over two-small elements and one full-sized element.}
\label{fig:SRD-1D}
\end{figure}

Unlike cell merging and linking, a given element may contribute to multiple merge neighborhoods. We denote all merge neighborhoods to which an element $D^k$ contributes as the set 
\begin{equation}
C_k = \{j : k \in M_j\}.
\end{equation}
We use the notation $|C_k|$ to denote the cardinality of $C_k$, i.e., number of elements in $C_k$. The final solution on each element is the average of all projected solutions to which a given element contributed. This process is illustrated on a 1D mesh with two small elements and one full-sized element in Figure \ref{fig:SRD-1D}.

State redistribution has been shown to relax the CFL condition of cut meshes to be proportional to the CFL condition of the underlying background mesh. Importantly, state redistribution preserves polynomial order and conserves the average of the solution \cite{berger-stateRedistr}. For systems of partial differential equations, such as the wave equation, state redistribution is applied to each component individually. State redistribution can be represented as a linear operator, $\SRDc$, acting on a semi-discrete solution $\bm{u}$, and is slightly dissipative in nature, making it a good partner for an energy stable scheme. The main computational complexity in using state redistribution is in setup: determining which elements are in need of stabilization, building their merge neighborhoods, and constructing the necessary operators. For a static mesh the state redistribution operators need only be computed once at the start of a simulation as a preprocessing step and then applied to the solution are each time step.

\subsection{Proof of $L_2$ Stability of High-Order State Redistribution}\label{sec:proof} % ----------------------------

To prove the energy stability of our final DG formulation with state redistribution, we appeal to work by Nordstr\"om and Winters \cite{nordstrom-StableFiltering} which states that an energy stable DG scheme to which a contractive filter is applied remains energy stable. Note that in general, filters can disrupt the energy stability of their underlying method \cite{lunquist-StableFiltering}. We interpret the state redistribution operator as a filter, reducing our proof to showing that the state redistribution operator is contractive in the $L_2$ norm.

\begin{theorem}
    Let $\SRDc$ be the state redistribution operator and $u \in L_2(\Omega)$. Then $\SRDc$ is contractive,
    \begin{equation}
    \left\| \SRDc u \right\|_{L_2(\domain)} \le \left\| u \right\|_{L_2(\domain)}.
\end{equation}
\end{theorem}\label{thm:1}

Before proving state redistribution is contractive we first introduce some notation and recall two basic lemmas. 
% The weighted inner product and norm
\begin{definition}
    We define the weighted inner product $(\cdot, \cdot)_{M_k}$ over the elements of a merge neighborhood $M_k$,
    $$D^{M_k} = \{D^j : j \in M_k\} \subseteq \Omega ,$$
    for $u,v \in \mathbb{P}^N(D^{M_k})$ as
    \begin{equation}
        (u,v)_{M_k} = \sum_{j \in M_k}^{}\frac{1}{|C_j|} (u,v)_{L_2(D^j)}
    \end{equation}
    and denote its associated norm as
    \begin{equation}
        \|u\|_{M_k} = (u,u)_{M_k} = \sum_{j \in M_k}^{}\frac{1}{|C_j|} \|u\|_{L_2(D^j)}.
    \end{equation}
\end{definition}

\begin{lemma}
\textit{(Jensen's inequality)} For a set of numbers $\{x_i\}_{i=1}^n \in \realR$ the following inequality holds \cite{jensens}
    $$\left( \sum_{i=1}^n x_i \right)^2 \le n \sum_{i=1}^n x_i^2.$$\label{lemma:Jensens}
\end{lemma}

\begin{lemma}
    Let $(\cdot, \cdot)$ denote some inner product on a domain $D$ with inner product norm $\|\cdot\|$. For a  given polynomial degree $N$, let $\Pi:\mathbb{P}^N(D)\rightarrow \mathbb{P}^N(D)$ be a projection operator satisfying
    $$(\Pi u,v) = (u,v), \quad \forall v \in \mathbb{P}^N(D).$$
    Then the projection operator $\Pi$ is contractive:
    $$\left\| \Pi u \right\| \le \|u\|.$$\label{lemma:contractivity}
\end{lemma}

The proofs of Lemmas \ref{lemma:Jensens} and \ref{lemma:contractivity} are straightforward using standard tools such as Jensen's and the Cauchy-Schwarz inequalities. We now present the proof of Theorem \ref{thm:1}.

\begin{proof}
    Let $(\cdot, \cdot)_{L_2(B)}$ denote the $L_2$ inner product on given subdomain $B \subseteq \Omega$. For a piecewise polynomial $u$ the action of merging the element-wise solutions from all elements in a merge neighborhood $M_k$ can be expressed as a volume-weighted $L_2$ projection $\Pi_k$ over $M_k$. Since state redistribution conserves the moments of the solution, $\Pi_k$ satisfies
    \begin{equation}
        (\Pi_k u, v)_{M_k} = (u,v)_{M_k}, \quad \forall v \in \mathbb{P}^N(D^{M_k}).
    \end{equation}
    The state redistribution operator restricted to $D^k$ can then be expressed as
    \begin{equation}
        \SRDc u|_{D^k} = \frac{1}{|C_k|}\sum_{j \in C_k} \Pi_j u.
    \end{equation}
    
    Using this notation the application of state redistribution over the entire domain can be expressed as
    \begin{align}
        \|\SRDc u\|_{L_2(\Omega)}^2 &= \sum_{k=1}^{N_h} \|\SRDc u|_{D^k}\|_{L_2(D^k)}^2 \\
        &= \sum_{k=1}^{N_h} \int_{D^k} \LRp{\frac{1}{|C_k|}\sum_{j \in C_k} \Pi_j u}^2 \\[6pt]
        &= \sum_{k=1}^{N_h} \frac{1}{|C_k|^2} \int_{D^k} \LRp{~\sum_{j \in C_k} \Pi_j u}^2.\label{eq:26}
    \end{align}
    Applying Lemma \ref{lemma:Jensens} to \eqref{eq:26}, we have that
    \begin{align}
        \sum_{k=1}^{N_h} \frac{1}{|C_k|^2} \int_{D^k} \LRp{~\sum_{j \in C_k} \Pi_j u}^2
        \le \sum_{k=1}^{N_h} \frac{1}{|C_k|}\sum_{j \in C_k} \|\Pi_j u\|_{L_2(D^k)}^2,\label{eqn:contrProof-ineq1}
    \end{align}
	where $N_h$ is the number of elements in the mesh.

    We can express nested sums over all $k$ and $j \in C_k$ equivalently as a nested sum over all $k$ and $j \in M_k$ and vice versa. The first of these corresponds to a sum over the final, merged solutions while the second corresponds to a sum of contributed terms from each element. Applying this change and the definition of the $\|\cdot \|_{M_k}$ norm yields
    \begin{align}
        \sum_{k=1}^{N_h} \frac{1}{|C_k|}\sum_{j \in C_k} \|\Pi_j u\|_{L_2(D^k)}^2 &= \sum_{k=1}^{N_h} \sum_{j \in M_k} \frac{1}{|C_j|}\|\Pi_k u\|_{L_2(D^j)}^2 \\
        &= \sum_{k=1}^{N_h} \|\Pi_k u\|_{M_k}^2.\label{eqn:weightedNorm}
    \end{align}
    We can then apply Lemma \ref{lemma:contractivity} to the RHS
    \begin{equation}
        \sum_{k=1}^{N_h} \|\Pi_k u\|_{M_k}^2 \le \sum_{M_k} \| u\|_{M_k}^2\label{eqn:contrProof-ineq2}
    \end{equation}
    before returning the sum of the weighted norm over all $k$ and $j \in M_k$ to sums of the $L_2$ norm over all $k$ and $j \in C_k$:
    \begin{equation}
        \sum_{k=1}^{N_h} \| u\|_{M_k}^2 = \sum_{k=1}^{N_h} \sum_{j \in M_k} \frac{1}{|C_j|} \|u\|^2_{L_2(D^j)} = \sum_{k=1}^{N_h} \frac{1}{|C_k|}\sum_{j \in C_k} \| u\|_{L_2(D^k)}^2.\label{eqn:4.36}
    \end{equation}
    Since the summand $\|u\|_{L_2(D^k)}$ does not depend on $j$ we have that
    \begin{equation}
        \sum_{j \in C_k} \| u\|_{L_2(D^k)}^2 =  |C_k| \| u\|_{L_2(D^k)}^2.
    \end{equation}
    Substituting this into Equation \eqref{eqn:4.36} gives 
    \begin{equation}
        \sum_{k=1}^{N_h} \frac{1}{|C_k|}\sum_{j \in C_k} \| u\|_{L_2(D^k)}^2 = \sum_{k=1}^{N_h} \| u\|_{L_2(D^k)}^2 = \|u\|_{L_2(\Omega)}^2
    \end{equation}
    which via the inequalities in \eqref{eqn:contrProof-ineq1} and \eqref{eqn:contrProof-ineq2} yield the desired result:
    $$\|\SRDc u\|_{L_2(\Omega)}^2 \le \|u\|_{L_2(\Omega)}^2.$$
\end{proof}

\section{Numerical Experiments}\label{sec:results} % ====================================================================

\subsection{Manufactured Solution}
To test the correctness of our solver we implement a manufactured solution for the wave equation
\begin{align}
\pd{p}{t}{} + \Div \bm{u} &= f(x,y,t)\label{eqn:mms-p} \\[10pt]
\pd{\bm{u}}{t}{} + \Grad p &= 0.\label{eqn:mms-u}
\end{align}
For the manufactured solution, we prescribe the pressure to be
\begin{equation}
p = \cos(2\pi t)\sin(\pi x)\sin(\pi y).
\end{equation}
Substituting the prescribed pressure into Equation \eqref{eqn:mms-u} and taking $\bm{u}(x,y,t=0) = 0$ yields
\begin{equation}
\bm{u} = \int_0^t \Grad p~ \text{d}\tau =  -\frac{1}{2}\begin{bmatrix}
~\sin(2 \pi t)\cos(\pi x)\sin(\pi y)~ \\[3pt] \sin(2\pi t)\sin(\pi x)\cos(\pi y)
\end{bmatrix}
\end{equation}
and
\begin{equation}
f(x,y,t) = -\pi \sin(2\pi t) \sin(\pi x)\sin(\pi y).
\end{equation}

We simulate this PDE (partial differential equation) on the Cartesian domain $[-1,1]^2$ from which we cut out a circle $C$ of radius $R=0.3$ centered at $(x_0, y_0)=(-0.5, 0)$. We use the manufactured solution to prescribe initial and boundary conditions. To assess our solver we performed a $h$-convergence study using degree $N=1, 2, 3,$ and $4$ polynomials and simulated the manufactured solution using $h = \Delta x = \Delta y = $ \nicefrac{1}{2}, \nicefrac{1}{4}, \nicefrac{1}{8}, \nicefrac{1}{16}. 

\begin{figure}[H]
\centering
.\includegraphics[width=0.75\textwidth]{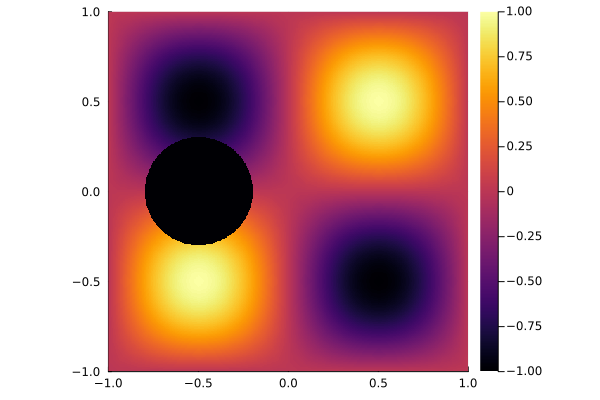}
\caption{Pressure field for the manufactured solution at $t=0$ showing the embedded object.}
\label{fig:MMS-IC}
\end{figure}

\begin{figure}[H]
\centering
\includegraphics[width=\textwidth]{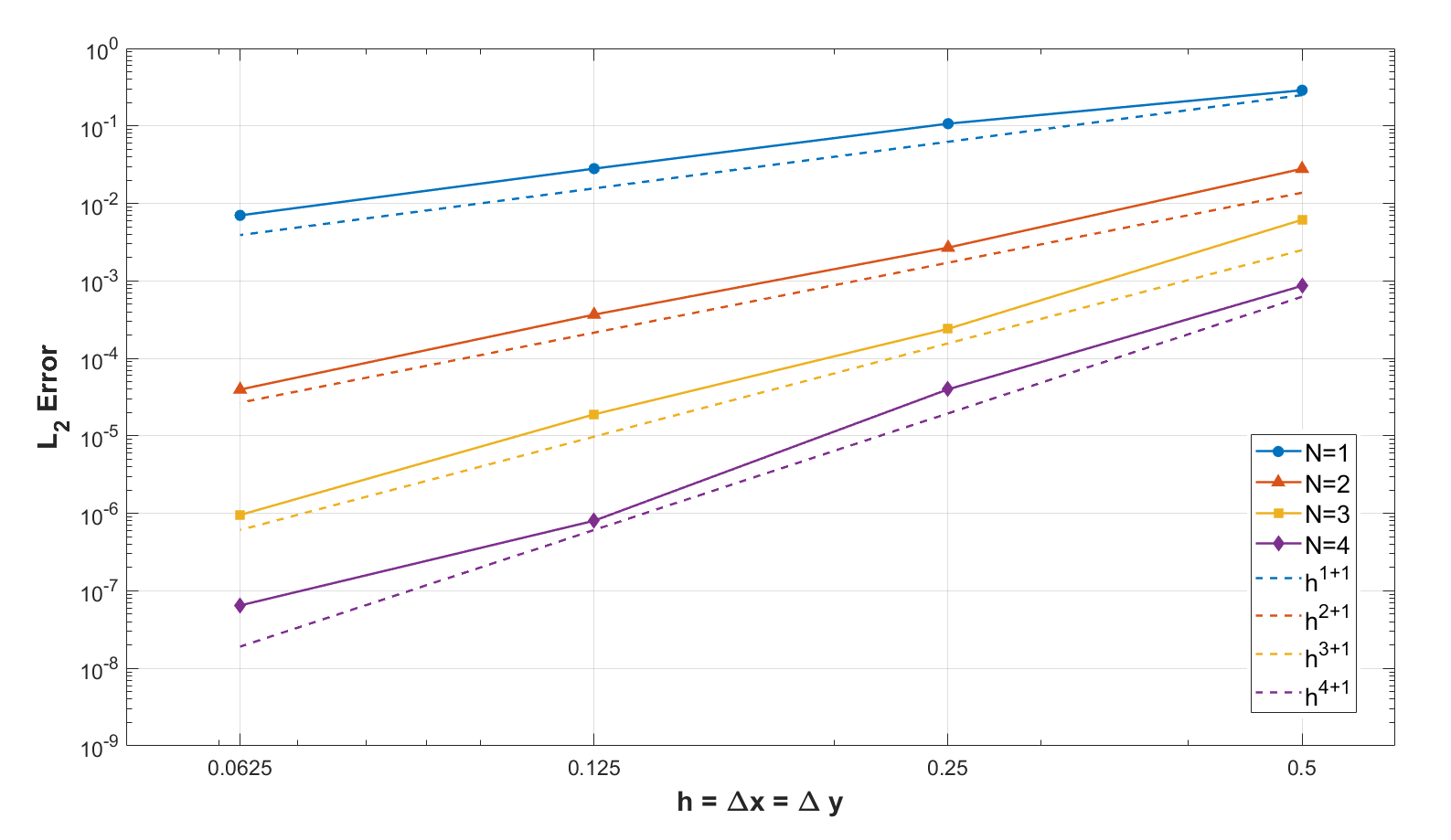}
\caption{$L_2$ error for the manufactured solution at various grid sizes, $h$, and polynomial degree, $N$, at time $t=1.3$.}
\label{fig:MMS-convergence}
\end{figure}

Figure \ref{fig:MMS-IC} shows the pressure solution at $t=0$. Note the velocity field is zero at $t=0$. For a degree $N$ solution, we expect $h^{N+1}$ convergence. Figure \ref{fig:MMS-convergence} shows the $L_2$ error convergence results with $h^{N+1}$ rates for comparison.

\subsection{Eigenvalues of RHS Operator}
In this experiment we seek to highlight the influence of state redistribution on CFL condition in the presence of a small cell by examining the eigenvalues of the RHS operator with the eigenvalues of the combined RHS-state redistribution operator. 

\medskip
After the applying the DG discretization to the wave equation, we are left with the global system of ordinary differential equations
\begin{equation}
\pd{\bm{U}}{t}{} = \bm{M}^{-1}\bm{QU} = \RHSd\bm{U}\label{eqn:RHS-op}
\end{equation}
where $\uSysC$ is the global vector of unknowns, $\bm{M}$ is the global mass matrix, $\bm{Q}$ is the discrete action of the spatial derivatives, and $\RHSd = \bm{M}^{-1}\bm{Q}$ is the overall action of the DG discretization prior to applying any time integration. 

State redistribution can be seen as a linear operator, which we will denote as $\SRDc$, that acts on the current solution vector. Adding state redistribution to our scheme modifies Equation \eqref{eqn:RHS-op} by applying state redistribution to the solution vector prior to applying the discretization operator: 
\begin{equation}
\pd{\bm{U}}{t}{} = \RHSc\SRDc\bm{U}\label{eqn:RHS-SRD-op}.
\end{equation}

\begin{figure}[H]
\centering
\includegraphics[width=\textwidth]{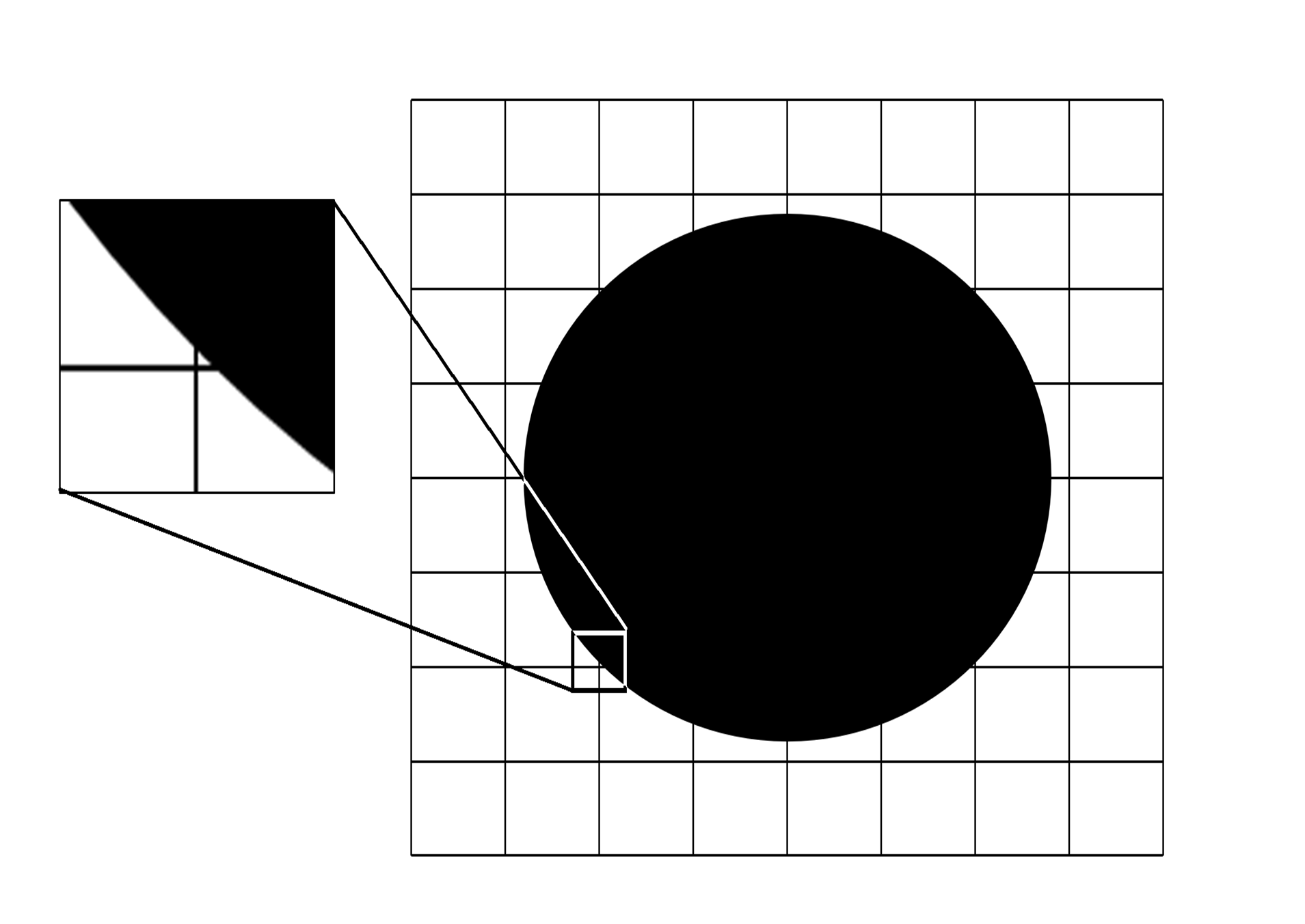}
\caption{The mesh and embedded object for the eigenvalue experiments, with one of smallest cells highlighted.}
\label{fig:eigval-mesh}
\end{figure}

We wish to compare the eigenvalues of $\RHSc$ with the eigenvalues of $\RHSc\SRDc$. To perform this comparison we consider a circle of radius $R=0.699$ embedded at the origin of a $8 \times 8$ Cartesian mesh on the domain $[-1,1]^2$. This setup yields a volume ratio of 940 between a full Cartesian cell and the smallest cut cell and a length ratio of 21.7. Note that we expect state redistribution to improve the CFL condition/maximum eigenvalue magnitude by a factor similar to the length ratio between a full cell and the smallest cut cell. For our solution space we use degree 4 polynomials. The mesh for the discretization is shown in Figure \ref{fig:eigval-mesh}.

Of interest in these experiments is the sign of the imaginary component of the eigenvalues and the maximum eigenvalue magnitude. The imaginary component of the eigenvalues reflects the energy conservation status of the scheme: positive imaginary component reflects increasing energy, zero imaginary component energy conservation, and negative imaginary component energy dissipation. The maximum eigenvalue magnitude meanwhile dictates the CFL condition of the scheme; the smaller the magnitude, the more gracious the CFL condition.

\begin{figure}[!h]
\centering
\begin{subfigure}[b]{0.75\textwidth}
        \centering
        \includegraphics[width=\textwidth]{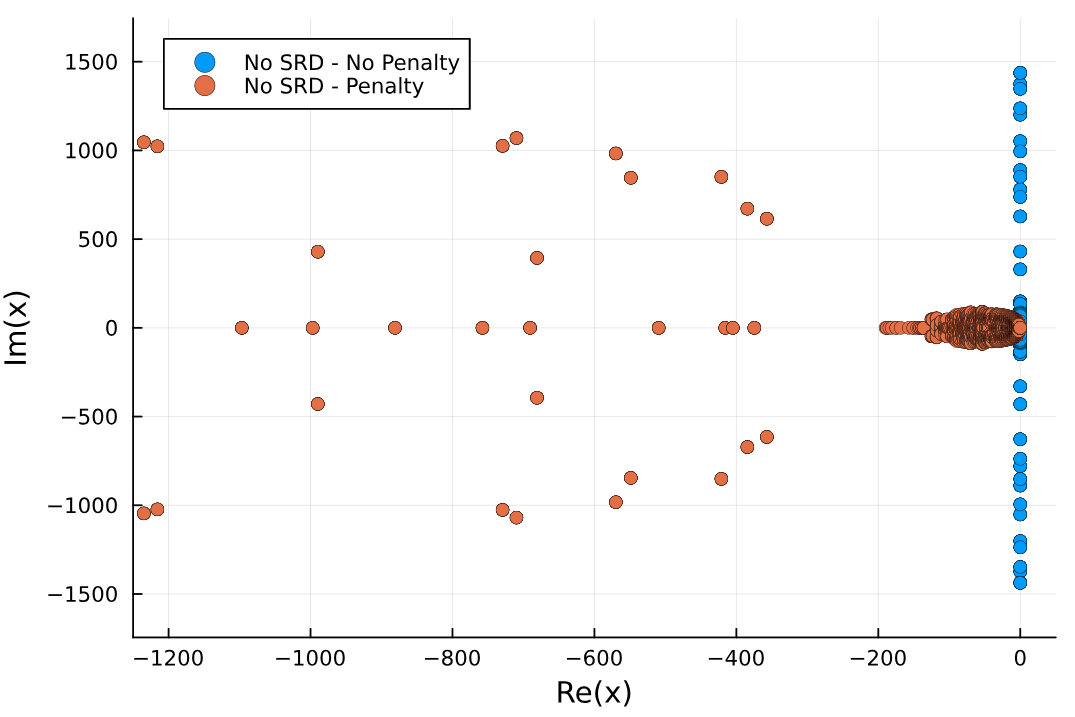}
        \caption{Eigenvalues without state redistribution.}\label{fig:eigvals-noSRD}
   \end{subfigure}

\begin{subfigure}[b]{0.75\textwidth}
        \centering
        \includegraphics[width=\textwidth]{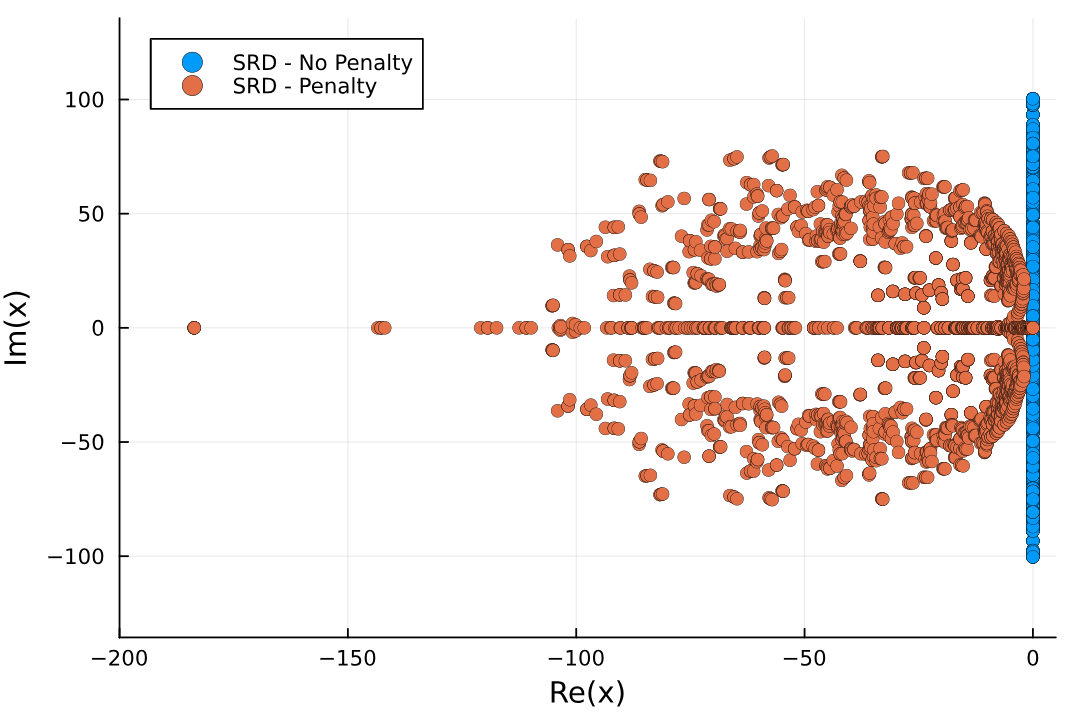}
        \caption{Eigenvalues with state redistribution.}\label{fig:eigvals-SRD}
   \end{subfigure}
\caption{Eigenvalues of the discretization operator with and without penalty terms and state redistribution. With penalty terms, the maximum magnitudes of the eignevalues' real components were 183 with state redistribution and 2157 without state redistribution. Without penalty terms, the maximum magnitudes of the real component of the eigenvalues were $6.07 \times 10^{-8}$ with state redistribution and $1.24 \times 10^{-13}$ without state redistribution. }
\label{fig:eigvals}
\end{figure}

Figure \ref{fig:eigvals-noSRD} shows the eigenvalues of the RHS operator without state redistribution for the original skew-symmetric energy conservative formulation and the penalized energy stable formulation.  With neither state redistribution nor penalty terms, the eigenvalues of the RHS operator fall solely on the real axis. When penalty terms are added, the eigenvalues' imaginary components are pushed into the left-half plane but maintain a large maximum magnitude. These results confirm the respective energy conservation and stability of these two schemes.

Figure \ref{fig:eigvals-SRD} shows the eigenvalues of the combined RHS-state redistribution operator with and without penalty terms. In both instances state redistribution is able to shrink the maximum eigenvalue magnitude and thus improve the CFL condition of the scheme. Table \ref{tab:eigvals} summarizes the changes in maximum eigenvalue magnitude both with and without state redistribution.

\begin{table}[!h]
\begin{center}
\caption{Maximum eigenvalues for the RHS operator with and without state redistribution and penalty terms.}
\begin{tabular}{lcc}
& No Penalty & Penalty \\ \hline
No State Redistribution & 1438 & 2157 \\
State Redistribution & 100 & 183\\
\hline
Ratio & 14.33 & 11.74 \\
\hline\hline\\[-9pt]
For reference:& Volume ratio=940 & Length ratio=21.7
\end{tabular}\label{tab:eigvals}
\end{center}
\end{table}

\subsection{Pacman Benchmark}
Next we compare our solver's solution with the ``Pacman" benchmark given in \cite{zeigelwanger-Pacman} to test the accuracy of our scheme. This benchmark provides an analytic solution to the acoustic wave equation for multiple wave sources scattering around a circle with an angular section removed, an object that resembles the video game character Pacman. This benchmark is of particular interest as it features a non-trivial embedded object and wave field while supplying an analytic solution to test simulations against. Here, we will only consider the case of a incoming plane wave originating from outside the pacman ``mouth".

For our simulation, we take the circle radius to be $R=1$ from which we remove the angular segment from $\theta = -\pi/6$ to $\theta = \pi / 6$. We simulate the resulting pacman object on the Cartesian domain $x \in [-3.3, 3]$, $y \in [-3,3]$ and a single $33 \times 33$ mesh. The offset in the lower $x$ boundary is to prevent split cut elements, where a single background Cartesian element is cut into multiple cut elements \cite{tao-tunneling}, which our code does not currently support.

As given in \cite{zeigelwanger-Pacman}, the analytic solution to the wave field is given in polar coordinates by infinite sums of either Hankel or Bessel functions depending on whether the point of interest is outside ($p^{(I)}, v_\theta^{(I)}, v_r^{(I)}$) or inside of the pacman mouth ($p^{(II)}, v_\theta^{(II)}, v_r^{(II)}$). The solution in the region outside the pacman mouth is given by:
\begin{align}
p^{(I)}(r,\theta) &= \sum_{n=0}^{\infty} \left(a_n^A \sin(n\theta) + a_n^S\cos(n\theta) \right)H_n^{(2)}(kr),\\[3pt]
v_r^{(I)}(r,\theta) &= \frac{\text{i}}{Z_0}\sum_{n=0}^{\infty} \left(a_n^A \sin(n\theta) + a_n^S\cos(n\theta) \right)H_n^{'(2)}(kr),\\[3pt]
v_\theta^{(I)}(r,\theta)& = \frac{\text{i}}{krZ_0}\sum_{n=0}^{\infty} n\left(a_n^A \cos(n\theta) - a_n^S\sin(n\theta) \right) H_n^{(2)}(kr).
\end{align}
Similarly the solution in the region inside the pacman mouth is given by:
\begin{align}
p^{(II)}(r,\theta) &= \sum_{n=0}^{\infty} b^A_n J_{(n + \nicefrac{1}{2})N}(kr)\sin\left(\left(n + \frac{1}{2}\right)N\theta \right) + b^S_n J_{nN}(kr)\cos\left(nN\theta \right),\\[3pt]
v_r^{(II)}(r,\theta)  &= \frac{\text{i}}{Z_0}\sum_{n=0}^{\infty} b^A_n J'_{(n + \nicefrac{1}{2})N}(kr)\sin\left(\left(n + \frac{1}{2}\right)N\theta \right) + b^S_n J'_{nN}(kr)\cos\left(nN\theta \right),\\[3pt]
v_\theta^{(II)}(r,\theta)  &= \frac{\text{iN}}{krZ_0}\sum_{n=0}^{\infty}  \left(n + \nicefrac{1}{2}\right)b^A_n J_{(n + \nicefrac{1}{2})N}(kr)\cos\left(\left(n + \frac{1}{2}\right)N\theta \right) - n b^S_n J_{nN}(kr)\sin\left(nN\theta \right),
\end{align}
where $H^{(2)}_n$ is the $n^{th}$ Hankel function of the second kind, $J_n$ is the $n^{th}$ Bessel function, $N = 6$ is the wedge number for defining the pacman mouth angle, and the coefficients $a^A_n, a^S_n, b^A_n, b^S_n$ are calculated as described in \cite{zeigelwanger-Pacman}. To compare this solution to our simulation, we truncate this sum to 100 terms. For our simulation we take $\rho = c_0 = 1$ and $k = 9$ which then dictates $Z_0 = 1$, $f = 1.432$, and $\omega = 9$. The plane wave is introduced from the upper right corner of the domain at an angle of $\nicefrac{\pi}{4}$.

To introduce the plane wave in our simulation, we use the analytic solution as the boundary conditions on the Cartesian domain boundaries and zero-velocity boundary conditions on the pacman boundaries. The analytic solution at $t=0$ is used as the initial condition for our solver. We ran the simulation up to an end time of $t = 1$ (approximately 1.43 periods) for degree $N=1,2,3,4$ polynomials using the \texttt{Tsit5} adaptive time integration scheme \cite{tsitouras-Tsit5} provided in \texttt{OrdinaryDiffEq.jl} \cite{Julia-OrdinaryDiff} and an initial time step of $\Delta t= 10^{-4}$. 

\begin{figure}[H]
\centering
\begin{subfigure}[b]{0.49\textwidth}
        \centering
        \includegraphics[width=\textwidth]{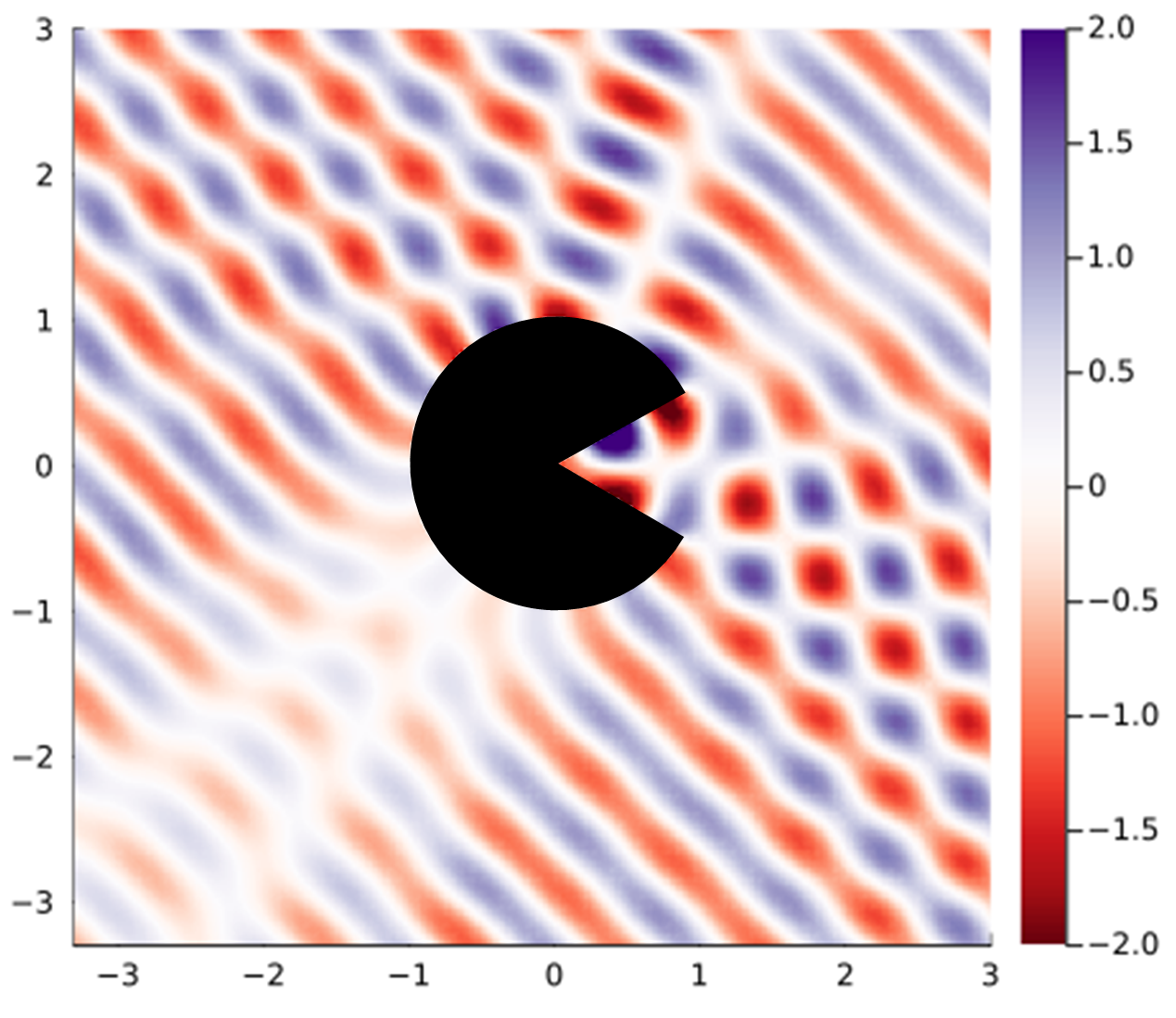}
        \caption{Analytic solution.}\label{fig:pacman-analytic}
\end{subfigure}
\begin{subfigure}[b]{0.49\textwidth}
        \centering
        \includegraphics[width=\textwidth]{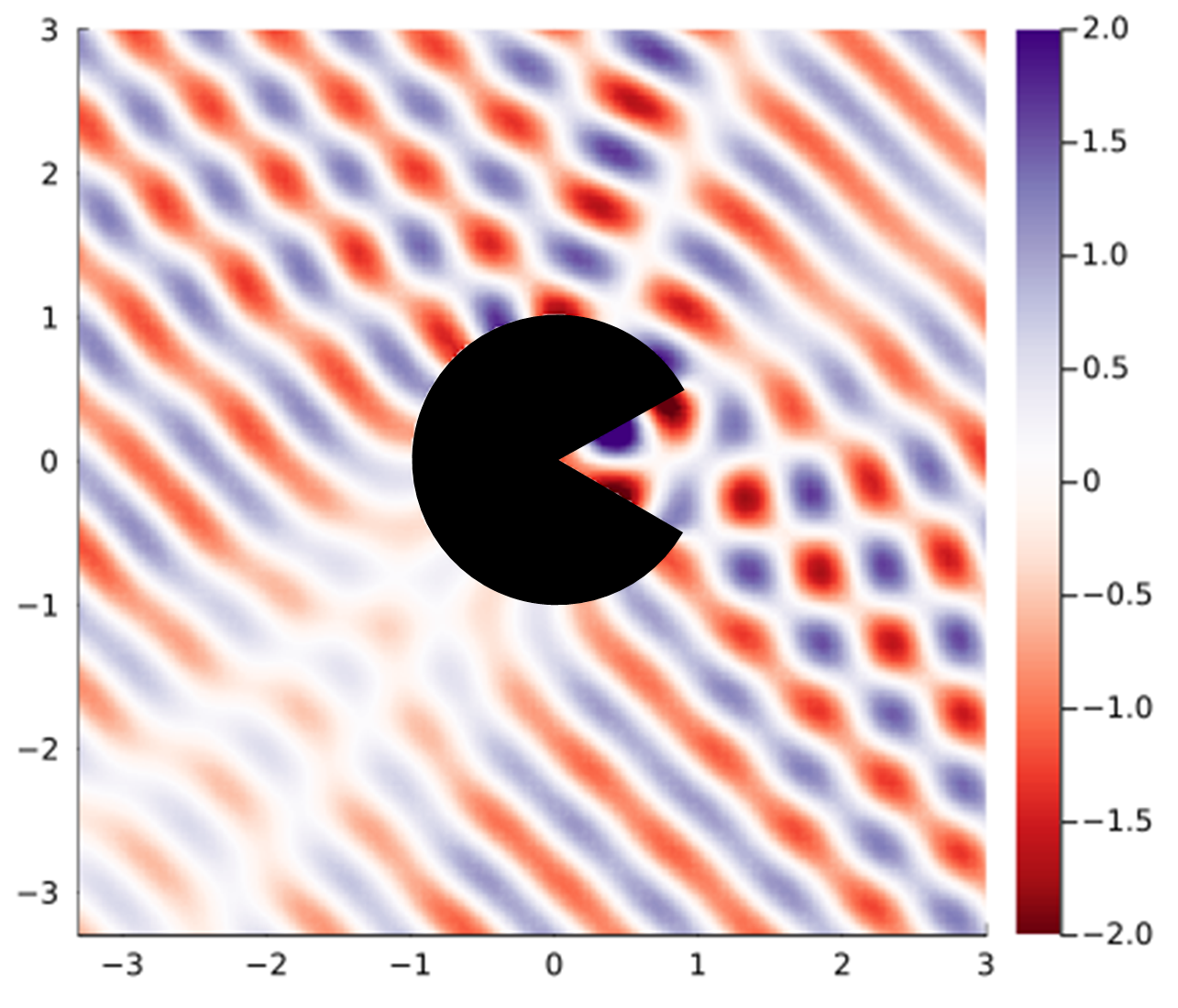}
        \caption{Solver solution for $N=4$.}\label{fig:pacman-analytic}
\end{subfigure}

\caption{Pressure field for the (a) analytic and (b) simulated solutions to the pacman wave scattering problem.}
\label{fig:pacman-solution}
\end{figure}

Figure \ref{fig:pacman-solution} shows the analytic and $N=4$ solution. While the two seem indistinguishable, there is in fact large localized error present in the solved solution. 

\begin{figure}[H]
\centering
\begin{subfigure}[b]{0.4\textwidth}
        \centering
        \includegraphics[width=\textwidth]{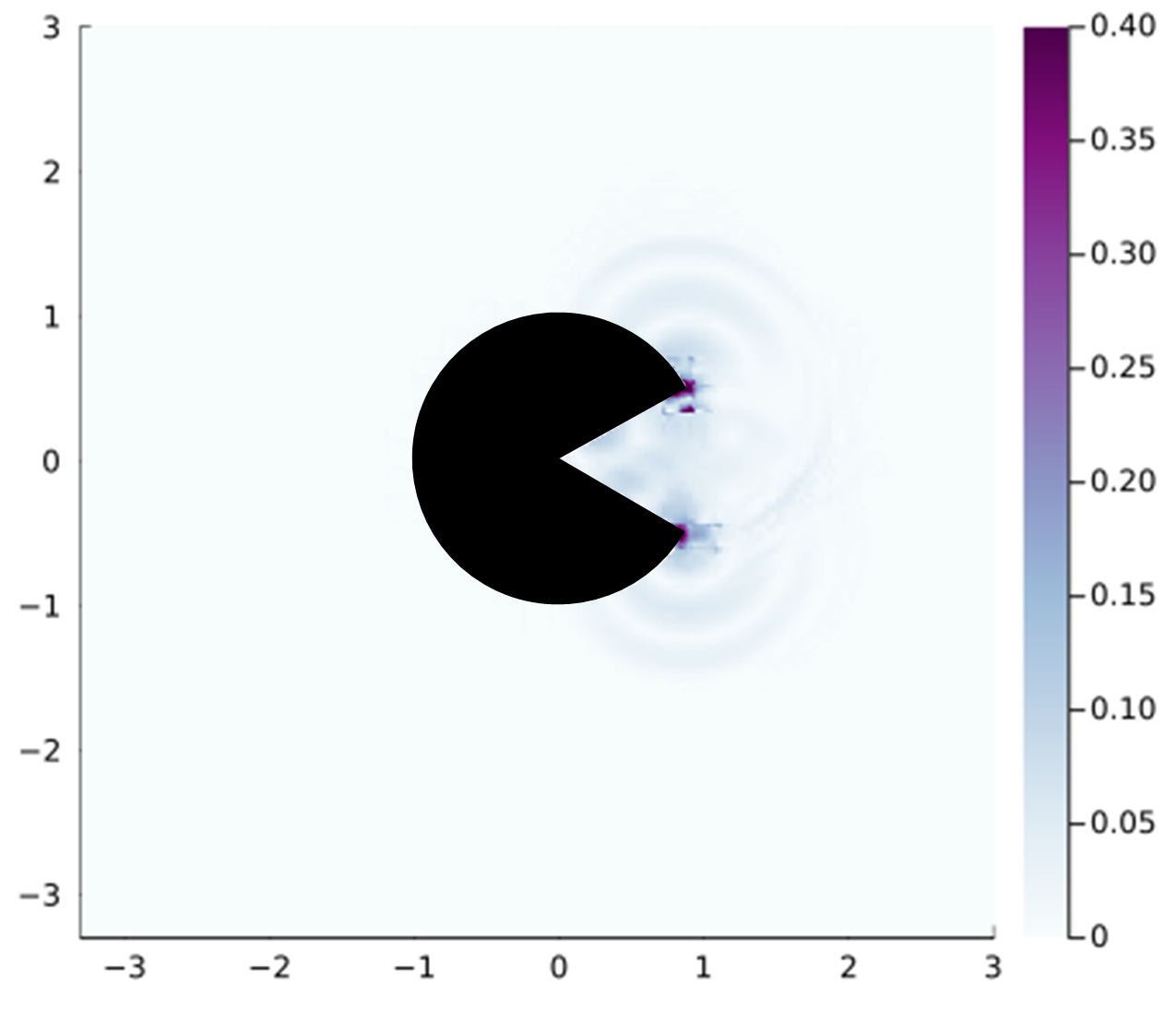}
        \caption{Threshold = $0.4$.}\label{fig:pacman-ptWsError-04}
\end{subfigure}
\hspace{5mm}
\begin{subfigure}[b]{0.4\textwidth}
        \centering
        \includegraphics[width=\textwidth]{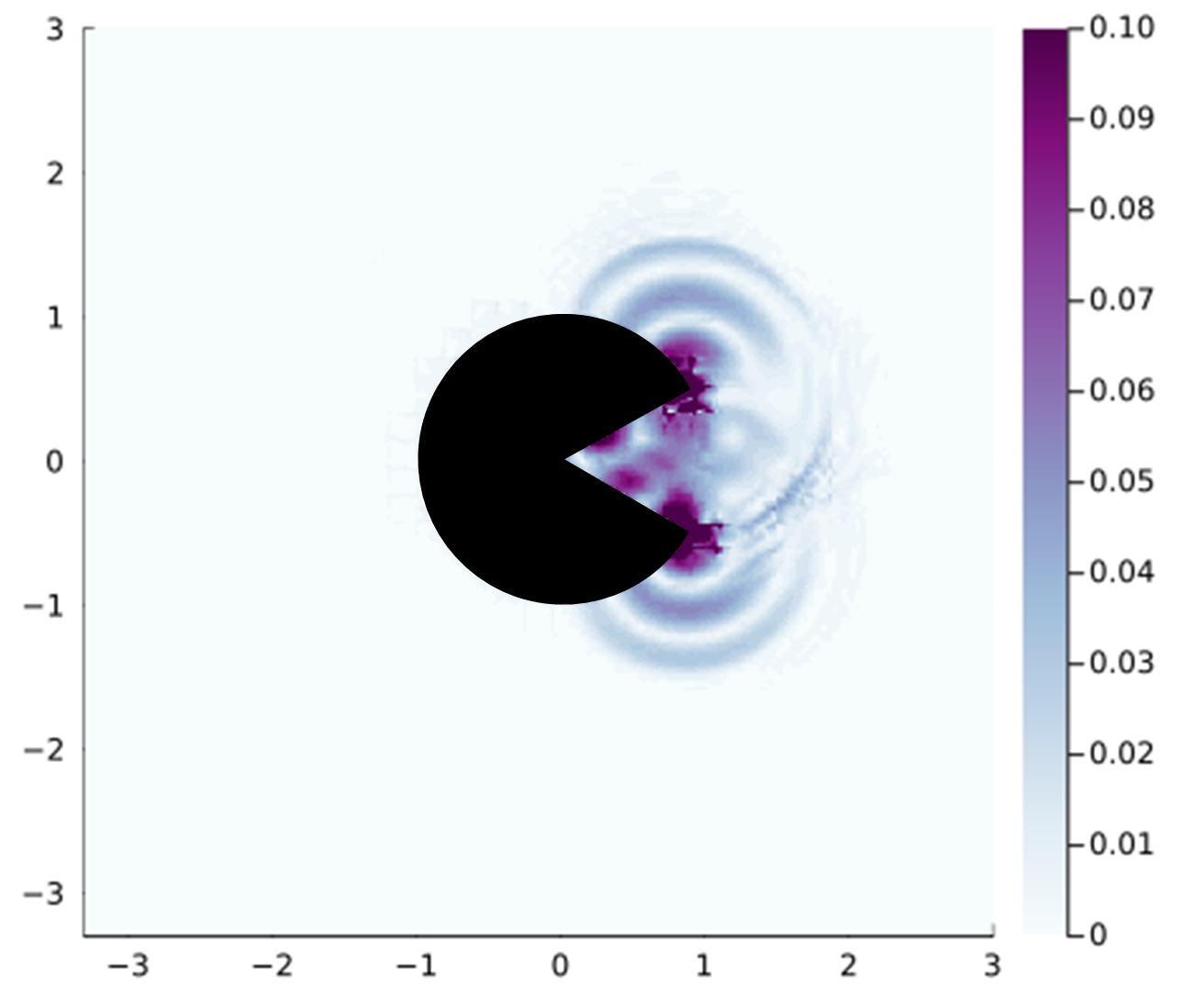}
        \caption{Threshold = $0.1$.}\label{fig:pacman-ptWsError-01}
\end{subfigure}

\vspace{5mm}
\begin{subfigure}[b]{0.6\textwidth}
        \centering
        \includegraphics[width=\textwidth]{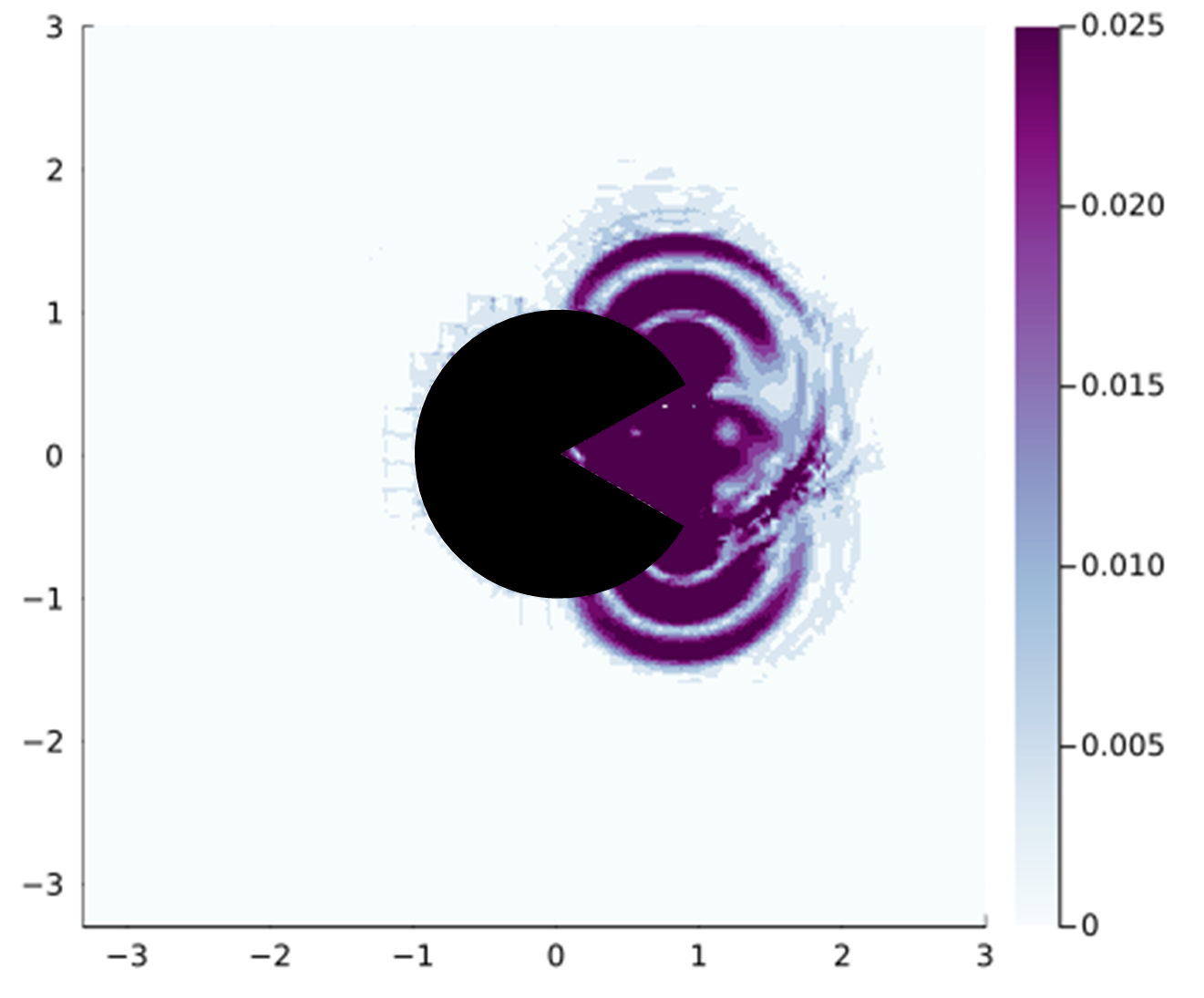}
        \caption{Threshold = $0.025$.}\label{fig:pacman-ptWsError-0025}
\end{subfigure}

\caption{Pointwise error at $t=1$ for the $N=4$ solution with the plotting scale clipped to various threshold values. }
\label{fig:pacman-error}
\end{figure}

As described in more detail in Appendix \ref{sec:Boundaries}, we represent embedded boundaries using explicit parameterizations. Due to this representation corners are represented exactly without smoothing. The sharpness of the convex corners is a source of localized error due to their lack of differentiability. Figure \ref{fig:pacman-error} shows the pointwise $\ell_2$ error in $(p, \bu)$ around the pacman object (clipped to various plotting thresholds). The actual maximum magnitude of the pointwise error is 0.96. As can be seen in Figure \ref{fig:pacman-ptWsError-04} (whose color scale saturates at 0.4), the largest error values are highly localized and propagate outward from the pacman mouth's corners (Figure \ref{fig:pacman-ptWsError-01}). Besides the error at the corners very low magnitude error can also be seen all around the embedded object in Figure \ref{fig:pacman-ptWsError-0025}. This low-level error around the object may be due to state redistribution and its dissipative nature.

\begin{figure}[H]
\centering
\begin{subfigure}[b]{0.44\textwidth}
        \centering
        \includegraphics[scale=0.65]{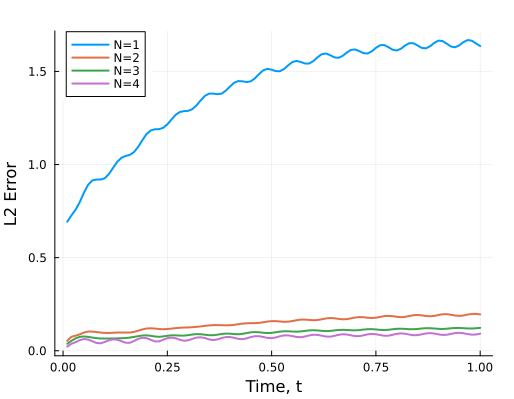}
        \caption{$L_2$ error.}\label{fig:pacman-L2}
\end{subfigure}
\hfill
\begin{subfigure}[b]{0.55\textwidth}
        \centering
        \includegraphics[scale=0.65]{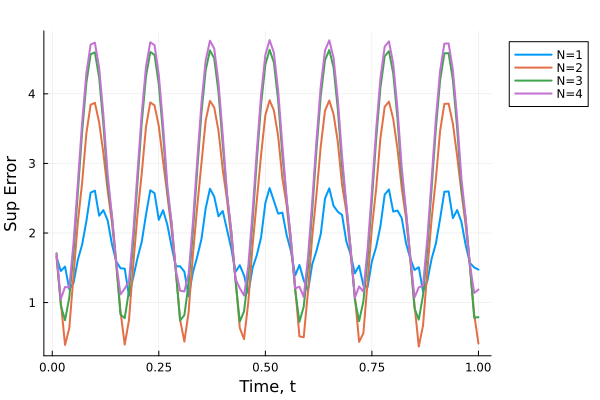}
        \caption{$L_\infty$ error.}\label{fig:pacman-Linf}
\end{subfigure}

\caption{$L_2$ and $L_\infty$ error versus time for the pacman benchmark. Note that despite the disturbance caused by the mouth corners the $L_\infty$ error remains bounded. }
\label{fig:pacman-errorVsTime}
\end{figure}

Figure \ref{fig:pacman-errorVsTime} shows the $L_2$ and $L_\infty$ error over time. Notably, the $L_\infty$ exhibits cyclic behavior; we theorize this is due to the corner effects of the pacman mouth. The $L_2$ errors also reflect this oscillatory behavior but still displays the expected long term behavior.

\subsection{Fish Simulation}
\begin{figure}[H]
\centering
\includegraphics[width=0.75\textwidth]{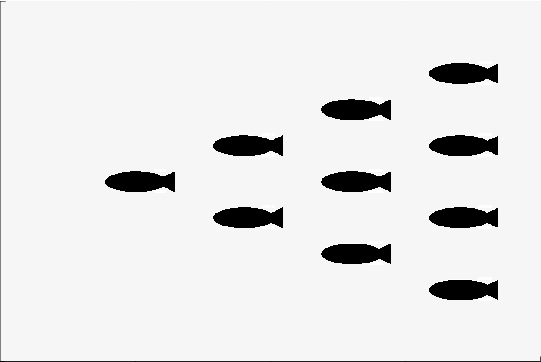}
\caption{Embedded objects for the fish simulation. Note the fish contain sharp concave and convex corners.}
\label{fig:fish-start}
\end{figure}

For our final numerical experiment we consider multiple embedded objects on a high-resolution mesh to test the versatility of our scheme. For this simulation, we embed 10 ``fish" arranged in a triangle pattern on the Cartesian domain $[-1, 1]^2$ over which we impose a $120 \times 120$ mesh. Figure \ref{fig:fish-start} shows the arrangement of fish in the Cartesian domain.

\begin{figure}[H]
\centering

\begin{subfigure}[b]{0.49\textwidth}
        \centering
        \includegraphics[width=0.95\textwidth]{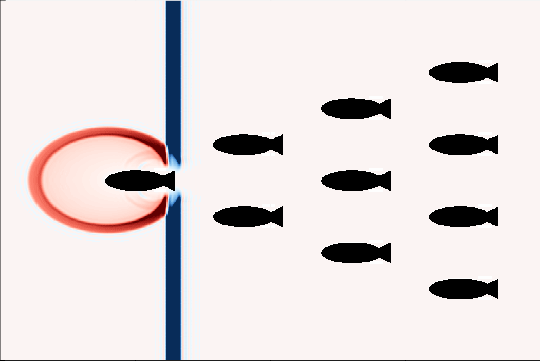}
        \caption{$t=0.12$.}\label{fig:fish-1}
\end{subfigure}
\begin{subfigure}[b]{0.49\textwidth}
        \centering
        \includegraphics[width=0.95\textwidth]{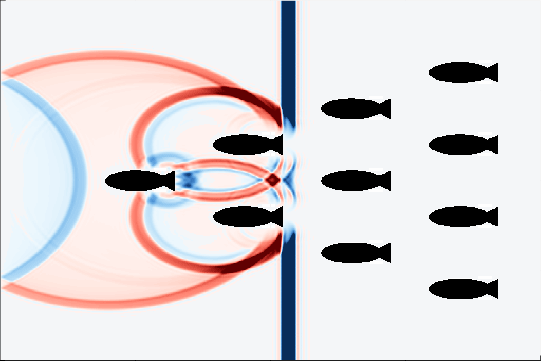}
        \caption{$t=0.19$.}\label{fig:fish-2}
\end{subfigure}

\vspace{5mm}
\begin{subfigure}[b]{0.49\textwidth}
        \centering
        \includegraphics[width=0.95\textwidth]{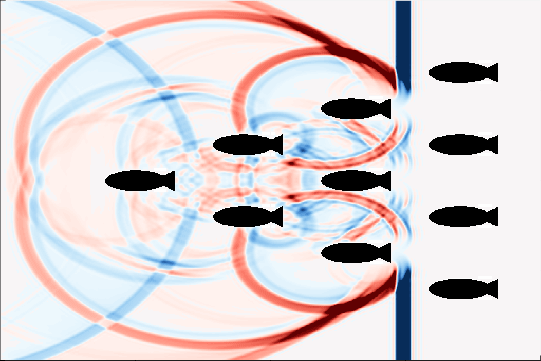}
        \caption{$t=0.26$.}\label{fig:fish-3}
\end{subfigure}
\begin{subfigure}[b]{0.49\textwidth}
        \centering
        \includegraphics[width=0.95\textwidth]{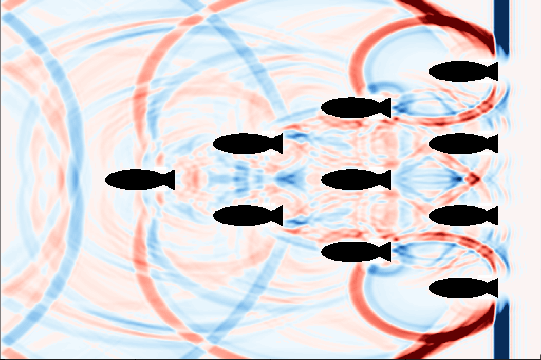}
        \caption{$t=0.32$.}\label{fig:fish-4}
\end{subfigure}

\caption{Pressure field at various times showing the reflection and scattering of the initial pressure wave off and through the fish.}
\label{fig:fish-snapshots}
\end{figure}

We impose zero-pressure boundary conditions at the right Cartesian boundary ($x = 1$) and fish. On the top and bottom Cartesian domain boundaries ($y = -1, 1$) we impose extrapolation boundary conditions. For times $ t \in (0, 0.05]$, the pressure on the left Cartesian boundary ($x = -1$) is set to $p = 2$ to generate a right-moving plane wave; for $t > 0.05$ zero-pressure conditions are enforced on this boundary instead. At time $t = 0$ both pressure and velocity fields are zero. We ran the simulation up to an end time of $t = 6$ using the \texttt{Tsit5} adaptive time integration scheme \cite{tsitouras-Tsit5} provided in \texttt{OrdinaryDiffEq.jl} \cite{Julia-OrdinaryDiff} and an initial time step of $\Delta t= 10^{-4}$. Figure \ref{fig:fish-snapshots} shows snapshots of the solution as the plane wave passes each column of fish at times $t = 0.12, 0.19, 0.26,$ and 0.32.

\section{Conclusions} %===============================================================================
We have proven that state redistribution can be added to an energy stable DG formulation without affecting the scheme's energy stability. We implemented a high-order energy stable DG solver with state redistribution for the acoustic wave equation on 2D cut meshes. To accommodate the loss of polynomially exact quadrature rules on cut elements we use a skew-symmetric DG formulation which is energy conservative under arbitrary quadrature and penalize the jumps along the interfaces of elements to add energy dissipation. Our code for constructing cut meshes and their associated DG and state redistribution operators are given in the registered Julia packages \texttt{PathIntersections.jl} and \texttt{StartUpDG.jl} \cite{Julia-PathIntersections, Julia-StartUpDG}. All code for our experiments is provided in this paper's accompanying GitHub repository \cite{code-repo}. Our codes represent embedded boundaries using explicit parameterizations which allows exact computation of normal vectors on the embedded boundaries and exact representation of sharp features. As shown in our numerical results, the presence of sharp geometry features can be a source of error due to their lack of differentiability. We have numerically verified our solver's expected convergence rates, energy stability, and relaxation of the CFL condition via state redistribution.

\section{Acknowledgments}
Christina Taylor and Jesse Chan gratefully acknowledges support from National Science Foundation under awards GRFP-1842494 and DMS-1943186. 
The authors also thank Dr.\ Marsha Berger and Dr.\ Andrew Giuliani for helpful discussions on state redistribution and constructing cut meshes.

\section{Appendix: Implementation Details}\label{sec:appendices} %==========================================================

\subsection{Embedded Boundary Representation}\label{sec:Boundaries}
To represent the embedded boundaries, we break with much of the Embedded Boundary method literature and use explicit parameterizations in place of the more common level-set functions. In the absence of a level-set/signed-distance function we determine which regions to include/exclude via a right-hand rule. The code to generate our cut mesh and its interpolation, quadrature, and state redistribution operators are provided in the respective Julia packages \texttt{PathIntersections.jl} and \texttt{StartUpDG.jl} \cite{Julia-PathIntersections, Julia-StartUpDG}. 

\begin{figure}[H]
\centering
\includegraphics[width=0.7\textwidth]{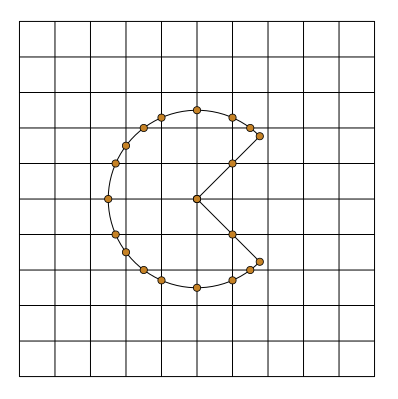}
\caption{Example of the mesh-curve intersections calculated by \texttt{PathIntersections.jl} for the pacman mesh. Note the stop points at the corners of the pacman mouth are included with the intersection points for use in defining element faces. }
\label{fig:PathIntersections}
\end{figure}

Given a parameterization of a curve $\curve_i:[0,1] \to \realR^2$ with parameter $s \in [0,1]$, our routines use a step-based algorithm to calculate mesh-curve intersections. This algorithm walks/``steps" along the curve from $s=0$ to $s=1$ sensing for intersections with the mesh and calculating the corresponding $s$ values as needed. This setup also allows users to define arbitrary ``stop points". Stop points can denote junctions in piecewise curves or points of interest on the curve to be used alongside intersection points to define the faces of cut elements. The intersections and stop points are then used to determine cut and excluded elements and build explicit piecewise representations of cut elements' boundaries. These routines are given in the Julia package \texttt{PathIntersections.jl}. More information on the routines can be found in the \texttt{PathIntersections.jl} documentation \cite{Julia-PathIntersections}. An example of the results of the intersection routine, including stop points, are shown in Figure \ref{fig:PathIntersections} for the pacman mesh.

The use of explicit parameterizations, while not extensible to 3D in our current implementation, comes with a number of benefits. Explicit parameterizations allow us to exactly represent sharp features and the curved boundaries of cut elements and exactly evaluate normal and tangent vectors via the use of automatic differentiation as provided in \texttt{ForwardDiff.jl} \cite{Julia-ForwardDiff}. These normal vectors are used elsewhere in our codes for constructing surface and volume quadratures.

\subsection{Construction of Merge Neighborhoods}
To determine when stabilization is needed for cut elements, we follow the convention in the original presentation of state redistribution in \cite{berger-stateRedistr} and stabilize cut elements whose volume is less than one-half the volume of the Cartesian elements. We can calculate the approximate volume of the cut elements using their volume quadrature rules applied to the constant function $f(x,y) = 1$.

To determine merge neighborhoods, we break with convention, which typically dictates merging with the nearest neighboring element in the direction normal to the embedded boundary. Instead, we use a greedy, volume-based approach. For a given cut element in need of merging, we compute its face neighbors. The neighbor that yields that largest increase in volume is added to the cut element's merge neighborhood. This is greedy process is repeated, now using the neighboring elements of the entire merge neighborhood, until the total volume of the cut neighborhood is at least threshold volume. Notably, the outcome of this process is affected by the order in which neighboring elements are found as multiple neighboring elements may have the same volume. In the case of such ties, the first element to be listed is used for merging.

There are other, more advanced ways of computing merge neighbors, such as incorporating the information on the outward normals of the cut boundary as previously mentioned. In many cases our greedy scheme, which will always favor merging with Cartesian elements over cut elements, resembles the merge neighborhoods that would be found via a normal-based scheme. This is especially true in the meshes we consider, where most cut elements have Cartesian or near Cartesian neighbors. In such meshes the greedy nature of our scheme makes it efficient with little impact on quality. For more complex meshes where cut elements have no Cartesian/near-Cartesian neighbors our scheme may result in oddly shaped/oblong merging neighborhoods which may impact the extent to which the CFL condition is relaxed.

\subsection{Construction of Volume Quadrature on Cut Elements}\label{sec:VolumeQuadr}
To construct volume quadrature rules on cut elements we use the divergence theorem to first compute target integrals via the surface quadrature. We then compute approximate Fekete points to use as quadrature nodes and solve for quadrature weights over those nodes using the target integral values. This approach is a form of moment fitting as described in \cite{davis-momentFitting}.

\subsubsection{Computing the Target Integrals}
On each cut element $D^k$ we seek a high-accuracy volume quadrature rule for integrating an arbitrary polynomial function $g(x,y) \in \mathbb{P}^{2N}(D^k)$.  To calculate volume quadrature rules we use an orthogonal 2D Legendre polynomial basis. Since these basis functions are defined on $[-1,1]^2$  we first map $(x,y) \in D^k$ to $(\hat{x},\hat{y}) \in \hat{D} = [-1,1]^2$. This mapping consists of shifting the centroid of the cut element, which we compute as the average of the face quadrature nodes, to the origin and then scaling by the maximum distance between the centroid and face nodes: 
\begin{align}
   \hat{x} &= \frac{x - \bar{x}_0}{\max_i \sqrt{(x_i - \bar{x}_0)^2, (y_i - \bar{y}_0)^2} }, \quad \bar{x}_0 = \frac{1}{n_q}\sum_{i=0}^{n_q} x_i, \\[3pt]
   \hat{y} &= \frac{y - \bar{y}_0}{\max_i \sqrt{(x_i - \bar{x}_0)^2, (y_i - \bar{y}_0)^2} }, \quad \bar{y}_0 = \frac{1}{n_q}\sum_{i=0}^{n_q} y_i,
\end{align}
where $(x_i, y_i)$ are the surface quadrature nodes on the faces of the element. We then define the functional $F$ that maps $g$ to a vector field $F(g)$ with $\div F(g) = g$ using
\begin{equation}
    F(g)(\hat{x},\hat{y}) = 
    \begin{bmatrix}
        \frac{1}{2} \int_{-1}^{\hat{x}} g(s,\hat{y}) \dint{s} \\
        \frac{1}{2} \int_{-1}^{\hat{y}} g(\hat{x},r) \dint{r} 
    \end{bmatrix}.
\end{equation}
The divergence theorem then states that:
\begin{equation}
    \int_{\hat{D}} g \dint{A} = \int_{\hat{D}} \div F(g) \dint{A} = \int_{\partial \hat{D}} (F(g) \cdot \bhat{n}) \dint{s}.\label{eq:divThm}
\end{equation}

By using explicit parameterizations of the embedded boundaries we have access to exact normal vectors which allows us to exactly map line quadrature rules to the edges of cut cells. However, given the embedded boundaries are not (in general) polynomially defined, the resulting surface quadrature rules are not polynomially exact, though they may be high accuracy. Via \eqref{eq:divThm} we use a highly accurate surface quadrature rule ($4(N+1)$ Gaussian quadrature nodes, where $N$ is the solution polynomial degree) to calculate target volume integral for each basis function, $\psi_i$,  to yield the vector of target integrals, $\bm{b}$
\begin{equation}
    \bm{b} = \begin{bmatrix}
        \int_{\hat{D}} \psi_0 \dint{A} \\
        \vdots \\
        \int_{\hat{D}} \psi_N \dint{A} 
    \end{bmatrix}
    = \begin{bmatrix}
        \int_{\partial \hat{D}} (F(\psi_0) \cdot \bhat{n}) \dint{s} \\
        \vdots \\
        \int_{\partial \hat{D}} (F(\psi_N) \cdot \bhat{n}) \dint{s}
    \end{bmatrix}.
\end{equation}

\subsubsection{Computing Volume Quadrature Nodes}
To compute quadrature nodes, we start by sampling a large number of equally spaced points in the cut element, $\{(x_k, y_k)\}_{k=1}^{N_s}$. We use these points to build a tall, rectangular Vandermonde matrix $\bm{V}$
\begin{equation}
    \bm{V} = \begin{bmatrix}
        \psi_0(x_1, y_1) & \cdots & \psi_N(x_1, y_1) \\
        \psi_0(x_2, y_2) & \cdots & \psi_N(x_2, y_2) \\
        \vdots & \ddots & \vdots \\
        ~\psi_0(x_{N_s-1}, y_{N_s-1}) & \cdots & \psi_N(x_{N_s-1}, y_{N_s-1})~ \\
        \psi_0(x_{N_s}, y_{N_s}) & \cdots & \psi_N(x_{N_s}, y_{N_s}) \\
    \end{bmatrix}.
\end{equation}
We apply QR factorization with pivoting to the matrix $\bm{V}^T$. Since our integrand is total degree $2N$, we take the $(2N+1)(2N+2)/2$ square matrix of pivot columns of $\bV^T$ to produce a reduced, square Vandermonde matrix $\bm{\tilde{V}}^T$. As defined in \cite{sommariva-feketeQR}, the sampling points corresponding to the pivot columns of the QR factorization are approximate Fekete points. Fekete points are interpolation points that ensure the reduced Vandermonde matrix $\bm{\tilde{V}}^T$ is well conditioned. With the vector of target integrals and the reduced Vandermonde matrix we can then solve the linear system
\begin{equation}
    \bm{\tilde{V}}^T \bm{w} = \bm{b}\label{eqn:Fekete}
\end{equation}
for $\bm{w}$, the vector of quadrature weights. 

While our process for generating volume quadrature rules on cut elements is very practical and utilizes robust routines like QR factorization it has two major drawbacks: it can return negative quadrature weights and, on extremely small elements, poorly chosen nodes. We discuss both of these issues next.

\subsubsection{Conditioning of the Volume Quadrature Weights on Cut Elements}

There are methods for generating purely positive quadrature weights, such as in \cite{saye-Quadrature, saye-QuadraturePolynomials}, but they are more complicated and computationally expensive to compute. Here, the properties of our formulation do not require positive quadrature weights so we accept the possibility of negative quadrature weights but acknowledge their impact on numerical round-off. 

\begin{figure}[H]
\centering
\begin{subfigure}[b]{0.4\textwidth}
        \centering
        \includegraphics[width=\textwidth]{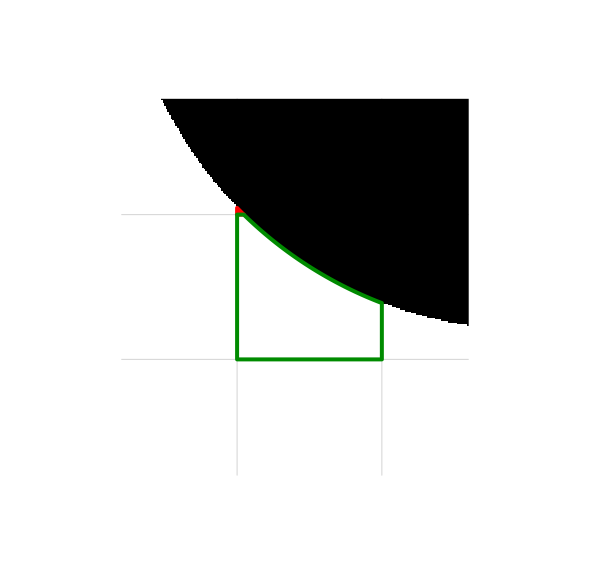}
        \caption{Circle mesh.}\label{fig:bestWorst-circle}
\end{subfigure}
\begin{subfigure}[b]{0.59\textwidth}
        \centering
        \includegraphics[width=\textwidth]{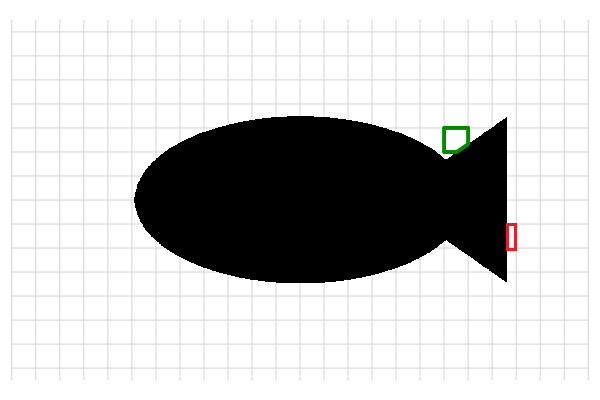}
        \caption{Fish mesh.}\label{fig:bestWorst-fish}
\end{subfigure}

\caption{Cut elements with the best (green) and worst (red) conditioned quadrature weights from the for $N=4$ schemes on the circle and fish meshes. In the case of the fish mesh, while the cut elements around each fish is analogous to the cut elements around other fish, analogous cut elements are slightly different. The best and worst elements actually occur on different fishes.}
\label{fig:bestWorst-quadr}
\end{figure}

\begin{figure}[H]
\centering
\includegraphics[width=0.7\textwidth]{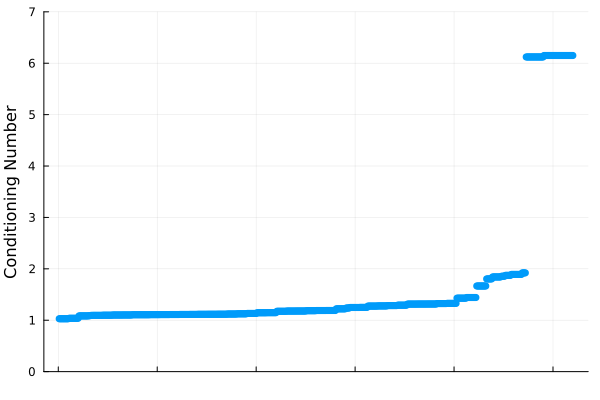}
\caption{Cut element quadrature weight conditioning numbers for the fish mesh for $N=4$ in order from best to worst. Of the 520 cut elements in this mesh 48 (9.2\%) have a conditioning number greater than 6. All other conditioning numbers are less than 2. }
\label{fig:quadCondNum-fish}
\end{figure}

In Table \ref{tab:volQuadr-cond} we list the conditioning of these quadrature weights for the circle mesh used in eigenvalue experiment and the fish mesh for polynomial degree $N=1,2,3,4$. These condition numbers are computed using the formula
\begin{equation}
\kappa_N = \LRp{\sum\limits_{i=1}^{m_N} |w_i|}\LRp{\sum\limits_{i=1}^{m_N} w_i}^{-1}
\end{equation}
where $w_i$ are the quadrature weights and $m_N$ is the number of quadrature points for the degree $N$ scheme. Note if all quadrature weights are non-negative the conditioning number will be 1. Figure \ref{fig:bestWorst-quadr} shows the cut elements with the best and worst conditioned quadrature weights for $N=4$ in each mesh while Figure \ref{fig:quadCondNum-fish} shows the spread of conditioning numbers for every cut element in the fish mesh.

\begin{table}[!h]
\centering
\caption{Best and worst cut element quadrature weight conditioning numbers for the circle and fish meshes.}\label{tab:volQuadr-cond}
\begin{tabular}{c|c|cc|cc}
 & & \multicolumn{2}{c}{\textbf{Circle Mesh}} & \multicolumn{2}{c}{\textbf{Fish Mesh}}\\
$N$ & $m_N$ & Best & Worst & Best & Worst \\
\hline
1 & 7  &  1      & 2.24163 & 1       & 2.17636 \\
2 & 16 & 1.19370 & 2.02327 & 1.03374 & 1.96366 \\
3 & 29 & 1.10240 & 2.63219 & 1.00875 & 2.64212 \\
4 & 46 & 1.08764 & 3.19112 & 1.02858 & 6.15015 \\
\hline\hline
\end{tabular}
\end{table}

\subsubsection{Poorly Chosen Quadrature Nodes}
The second, more problematic issue of poorly chosen quadrature nodes arises with extremely small elements. Our routine finds the sampling points closest to true Fekete points. As a result the quality of the final approximate Fekete points depends on the resolution of the sampling points. We generate sampling points via an evenly spaced grid of points on the background Cartesian element from which we remove points that are not in the cut element. If the cut element is extremely small, very few of the generated points will actually be in the cut element and thus represented in $\bm{V}^T$, from which the Fekete points are calculated. This issue is reflected in the rank of $\bm{V}^T$: if there where not a sufficient number of well-positioned points, $\bm{V}^T$ will be low-rank. 

To address this issue, we check the conditioning of $\bm{V}^T$ and resample points on a finer grid until $\bm{V}^T$ has sufficient rank. This resampling likely influences the clustering of conditioning number values seen in Figure \ref{fig:quadCondNum-fish}. However, even with this safe guard, the conditioning of the quadrature rules on extremely small cut cells may be poor. There are alternative methods for generating volume quadrature such as work by Saye in \cite{saye-Quadrature} and \cite{saye-QuadraturePolynomials}. However, for this proof of concept work we used this approach as it is relatively straightforward in comparison to these other methods.

\bibliographystyle{elsarticle-num}
\bibliography{bibliography}

\begin{thebibliography}{10}
\expandafter\ifx\csname url\endcsname\relax
  \def\url#1{\texttt{#1}}\fi
\expandafter\ifx\csname urlprefix\endcsname\relax\def\urlprefix{URL }\fi
\expandafter\ifx\csname href\endcsname\relax
  \def\href#1#2{#2} \def\path#1{#1}\fi

\bibitem{berger-stateRedistr}
M.~Berger, A.~Giuliani, A state redistribution algorithm for finite volume
  schemes on cut cell meshes, Journal of Computational Physics 428 (2021)
  109820.
\newblock \href {https://doi.org/10.1016/j.jcp.2020.109820}
  {\path{doi:10.1016/j.jcp.2020.109820}}.

\bibitem{dafermos-HyperbolicConservationLawsBook}
C.~M. Dafermos, Hyperbolic Conservation Laws in Continuum Physics, A Series of
  Comprehensive Studies in Mathematics, Springer-Verlag Berlin Heidelberg,
  2009.
\newblock \href {https://doi.org/10.1007/978-3-642-04048-1}
  {\path{doi:10.1007/978-3-642-04048-1}}.

\bibitem{warburton-skewSymm}
T.~Warburton, A low-storage curvilinear discontinuous {G}alerkin method for
  wave problems, SIAM Journal on Scientific Computing 35~(4) (2013)
  A1987--A2012.
\newblock \href {https://doi.org/10.1137/120899662}
  {\path{doi:10.1137/120899662}}.

\bibitem{wang-highVsLowOrder-HOinstability}
Z.~Wang, K.~Fidkowski, R.~Abgrall, F.~Bassi, D.~Caraeni, A.~Cary, H.~Deconinck,
  R.~Hartmann, K.~Hillewaert, H.~Huynh, N.~Kroll, G.~May, P.-O. Persson, B.~van
  Leer, M.~Visbal, High-order {CFD} methods: current status and perspective,
  International Journal for Numerical Methods in Fluids 72~(8) (2013) 811--845.
\newblock \href {https://doi.org/10.1002/fld.3767}
  {\path{doi:10.1002/fld.3767}}.

\bibitem{visbal-highOrder-unsteadyFlows}
M.~R. Visbal, D.~V. Gaitonde, High-order-accurate methods for complex unsteady
  subsonic flows, AIAA Journal 37~(10) (1999) 1231--1239.
\newblock \href {https://doi.org/10.2514/2.591} {\path{doi:10.2514/2.591}}.

\bibitem{ainsworth-highOrder-dispersionError}
M.~Ainsworth, Dispersive and dissipative behaviour of high order discontinuous
  {G}alerkin finite element methods, Journal of Computational Physics 198~(1)
  (2004) 106--130.
\newblock \href {https://doi.org/10.1016/j.jcp.2004.01.004}
  {\path{doi:10.1016/j.jcp.2004.01.004}}.

\bibitem{berger-cutStabilityAccuracy}
M.~Berger, Chapter 1 - {C}ut cells: Meshes and solvers, in: R.~Abgrall, C.-W.
  Shu (Eds.), Handbook of Numerical Methods for Hyperbolic Problems, Vol.~18 of
  Handbook of Numerical Analysis, Elsevier, 2017, pp. 1--22.
\newblock \href {https://doi.org/10.1016/bs.hna.2016.10.008}
  {\path{doi:10.1016/bs.hna.2016.10.008}}.

\bibitem{reed-originalCutCell}
W.~H. Reed, T.~R. Hill, Triangular mesh methods for the neutron transport
  equation, Tech. Rep. {LA-UR-73-479}, Los Alamos Scientific Lab., N. Mex.(USA)
  (1973).

\bibitem{quirk-EarlyCutCell}
J.~J. Quirk, An alternative to unstructured grids for computing gas dynamic
  flows around arbitrarily complex two-dimensional bodies, Computers \& Fluids
  23~(1) (1994) 125--142.
\newblock \href {https://doi.org/10.1016/0045-7930(94)90031-0}
  {\path{doi:10.1016/0045-7930(94)90031-0}}.

\bibitem{ingram-cutCellDevelopments}
D.~Ingram, D.~Causon, C.~Mingham, Developments in {C}artesian cut cell methods,
  Mathematics and Computers in Simulation 61~(3) (2003) 561--572, mODELLING
  2001 - Second IMACS Conference on Mathematical Modelling and Computational
  Methods in Mechanics, Physics, Biomechanics and Geodynamics.
\newblock \href {https://doi.org/10.1016/S0378-4754(02)00107-6}
  {\path{doi:10.1016/S0378-4754(02)00107-6}}.

\bibitem{chung-cellMerging}
M.-H. Chung, {C}artesian cut cell approach for simulating incompressible flows
  with rigid bodies of arbitrary shape, Computers \& Fluids 35~(6) (2006)
  607--623.
\newblock \href {https://doi.org/10.1016/j.compfluid.2005.04.005}
  {\path{doi:10.1016/j.compfluid.2005.04.005}}.

\bibitem{schott-cutFSIMonolithic}
B.~Schott, C.~Ager, W.~A. Wall, Monolithic cut finite element–based
  approaches for fluid-structure interaction, International Journal for
  Numerical Methods in Engineering 119~(8) (2019) 757--796.
\newblock \href {https://doi.org/10.1002/nme.6072}
  {\path{doi:10.1002/nme.6072}}.

\bibitem{udaykumar-phaseFront}
H.~Udaykumar, R.~Mittal, W.~Shyy, Computation of solid–liquid phase fronts in
  the sharp interface limit on fixed grids, Journal of Computational Physics
  153~(2) (1999) 535--574.
\newblock \href {https://doi.org/10.1006/jcph.1999.6294}
  {\path{doi:10.1006/jcph.1999.6294}}.

\bibitem{udaykumar-movingCutCell}
H.~Udaykumar, H.-C. Kan, W.~Shyy, R.~Tran-Son-Tay, Multiphase dynamics in
  arbitrary geometries on fixed {C}artesian grids, Journal of Computational
  Physics 137~(2) (1997) 366--405.
\newblock \href {https://doi.org/10.1006/jcph.1997.5805}
  {\path{doi:10.1006/jcph.1997.5805}}.

\bibitem{berre-rognes}
N.~Berre, M.~E. Rognes, A.~Massing, Cut finite element discretizations of
  cell-by-cell {EMI} electrophysiology models (2023).
\newblock \href {http://arxiv.org/abs/2306.03001} {\path{arXiv:2306.03001}}.

\bibitem{chan-skewSymm}
J.~Chan, Z.~Wang, A.~Modave, J.-F. Remacle, T.~Warburton, {GPU}-accelerated
  discontinuous {G}alerkin methods on hybrid meshes, Journal of Computational
  Physics 318 (2016) 142--168.
\newblock \href {https://doi.org/10.1016/j.jcp.2016.04.003}
  {\path{doi:10.1016/j.jcp.2016.04.003}}.

\bibitem{xie-ImplicitCut}
Z.~Xie, An implicit {C}artesian cut-cell method for incompressible viscous
  flows with complex geometries, Computer Methods in Applied Mechanics and
  Engineering 399 (2022) 115449.
\newblock \href {https://doi.org/10.1016/j.cma.2022.115449}
  {\path{doi:10.1016/j.cma.2022.115449}}.

\bibitem{may-explicitImplicit-original}
S.~May, M.~Berger, A mixed explicit implicit time stepping scheme for
  {C}artesian embedded boundary meshes, in: J.~Fuhrmann, M.~Ohlberger, C.~Rohde
  (Eds.), {Finite Volumes for Complex Applications VII-Methods and Theoretical
  Aspects}, Springer International Publishing, Cham, 2014, pp. 393--400.
\newblock \href {https://doi.org/10.1007/978-3-319-05684-5\_38}
  {\path{doi:10.1007/978-3-319-05684-5\_38}}.

\bibitem{may-explicitImplicit-accuracy}
S.~May, M.~Berger, F.~Laakmann, Accuracy considerations of mixed explicit
  implicit schemes for embedded boundary meshes, Proceedings in Applied
  Mathematics and Mechanics 19~(1) (2019) e201900411.
\newblock \href {https://doi.org/10.1002/pamm.201900411}
  {\path{doi:10.1002/pamm.201900411}}.

\bibitem{schoeder-localTimestepping}
S.~Schoeder, S.~Sticko, G.~Kreiss, M.~Kronbichler, High-order cut discontinuous
  {G}alerkin methods with local time stepping for acoustics, International
  Journal for Numerical Methods in Engineering 121~(13) (2020) 2979--3003.
\newblock \href {https://doi.org/10.1002/nme.6343}
  {\path{doi:10.1002/nme.6343}}.

\bibitem{sticko-stabilizedBilinearForm}
S.~Sticko, G.~Ludvigsson, G.~Kreiss, High-order cut finite elements for the
  elastic wave equation, Advances in Computational Mathematics 46~(3) (2020)
  45.
\newblock \href {https://doi.org/10.1007/s10444-020-09785-z}
  {\path{doi:10.1007/s10444-020-09785-z}}.

\bibitem{fu-cutDGPenaltyTerms}
P.~Fu, G.~Kreiss, High order cut discontinuous {G}alerkin methods for
  hyperbolic conservation laws in one space dimension, SIAM Journal on
  Scientific Computing 43~(4) (2021) A2404--A2424.
\newblock \href {https://doi.org/10.1137/20M1349060}
  {\path{doi:10.1137/20M1349060}}.

\bibitem{gurkan-aPriori-geomIndep}
C.~G\"{u}rkan, S.~Sticko, A.~Massing,
  \href{https://doi.org/10.1137/18M1206461}{Stabilized cut discontinuous
  {G}alerkin methods for advection-reaction problems}, SIAM Journal on
  Scientific Computing 42~(5) (2020) A2620--A2654.
\newblock \href {https://doi.org/10.1137/18M1206461}
  {\path{doi:10.1137/18M1206461}}.
\newline\urlprefix\url{https://doi.org/10.1137/18M1206461}

\bibitem{sticko-nistche-aPrioriError}
S.~Sticko, G.~Kreiss, A stabilized {N}itsche cut element method for the wave
  equation, Computer Methods in Applied Mechanics and Engineering 309 (2016)
  364--387.
\newblock \href {https://doi.org/10.1016/j.cma.2016.06.001}
  {\path{doi:10.1016/j.cma.2016.06.001}}.

\bibitem{dePrenter-noteOnNistche}
F.~{de Prenter}, C.~Lehrenfeld, A.~Massing, A note on the stability parameter
  in {N}itsche’s method for unfitted boundary value problems, Computers \&
  Mathematics with Applications 75~(12) (2018) 4322--4336.
\newblock \href {https://doi.org/10.1016/j.camwa.2018.03.032}
  {\path{doi:10.1016/j.camwa.2018.03.032}}.

\bibitem{engwer-DoD-original}
C.~Engwer, S.~May, A.~N\"{u}\ss{}ing, F.~Streitb\"{u}rger, A stabilized {DG}
  cut cell method for discretizing the linear transport equation, SIAM Journal
  on Scientific Computing 42~(6) (2020) A3677--A3703.
\newblock \href {https://doi.org/10.1137/19M1268318}
  {\path{doi:10.1137/19M1268318}}.

\bibitem{birke-DoD-linearHyperbolic}
G.~Birke, C.~Engwer, S.~May, F.~Streitbürger, {DoD S}tabilization of linear
  hyperbolic {PDE}s on general cut-cell meshes, Proceedings of Apllied
  Mathematics and Mechanics 23~(1) (2023) e202200198.
\newblock \href {https://doi.org/10.1002/pamm.202200198}
  {\path{doi:10.1002/pamm.202200198}}.

\bibitem{may-DoD-nonlinearHyperbolic}
S.~May, F.~Streitbürger, {DoD S}tabilization for non-linear hyperbolic
  conservation laws on cut cell meshes in one dimension, Applied Mathematics
  and Computation 419 (2022) 126854.
\newblock \href {https://doi.org/10.1016/j.amc.2021.126854}
  {\path{doi:10.1016/j.amc.2021.126854}}.

\bibitem{streitburger-DoD-monotonicity}
F.~Streitb{\"u}rger, C.~Engwer, S.~May, A.~N{\"u}{\ss}ing, Monotonicity
  considerations for stabilized {DG} cut cell schemes for the unsteady
  advection equation, in: F.~J. Vermolen, C.~Vuik (Eds.), Numerical Mathematics
  and Advanced Applications ENUMATH 2019, Springer International Publishing,
  Cham, 2021, pp. 929--937.
\newblock \href {https://doi.org/10.1007/978-3-030-55874-1\_92}
  {\path{doi:10.1007/978-3-030-55874-1\_92}}.

\bibitem{berger-hBox-original}
M.~Berger, R.~LeVeque, {C}artesian meshes and adaptive mesh refinement for
  hyperbolic partial differential equations, in: Proceedings of the Third
  International Conference on Hyperbolic Problems. Uppsala, Sweden. June 1990,
  1990.

\bibitem{helzel-hBox}
C.~Helzel, M.~J. Berger, R.~J. Leveque, A high-resolution rotated grid method
  for conservation laws with embedded geometries, SIAM Journal on Scientific
  Computing 26~(3) (2005) 785--809.
\newblock \href {https://doi.org/10.1137/S106482750343028X}
  {\path{doi:10.1137/S106482750343028X}}.

\bibitem{berger-hBox-cut}
M.~J. Berger, C.~Helzel, R.~J. LeVeque, H-box methods for the approximation of
  hyperbolic conservation laws on irregular grids, SIAM Journal on Numerical
  Analysis 41~(3) (2003) 893--918.
\newblock \href {https://doi.org/10.1137/S0036142902405394}
  {\path{doi:10.1137/S0036142902405394}}.

\bibitem{berger-simplifiedhBox}
M.~Berger, C.~Helzel, A simplified h-box method for embedded boundary grids,
  SIAM Journal on Scientific Computing 34~(2) (2012) A861--A888.
\newblock \href {https://doi.org/10.1137/110829398}
  {\path{doi:10.1137/110829398}}.

\bibitem{chern-fluxRedistrbution-early}
I.-L. Chern, P.~Colella, A conservative front tracking method for hyperbolic
  conservation laws, LLNL Rep. No. UCRL-97200, Lawrence Livermore National
  Laboratory 51 (1987) 83--110.

\bibitem{colella-fluxRedistribution}
P.~Colella, D.~T. Graves, B.~J. Keen, D.~Modiano, A {C}artesian grid embedded
  boundary method for hyperbolic conservation laws, Journal of Computational
  Physics 211~(1) (2006) 347--366.
\newblock \href {https://doi.org/10.1016/j.jcp.2005.05.026}
  {\path{doi:10.1016/j.jcp.2005.05.026}}.

\bibitem{berger-CellMerging}
M.~Berger, A note on the stability of cut cells and cell merging, Applied
  Numerical Mathematics 96 (2015) 180--186.
\newblock \href {https://doi.org/10.1016/j.apnum.2015.05.003}
  {\path{doi:10.1016/j.apnum.2015.05.003}}.

\bibitem{muralidharan-cellLinking}
B.~Muralidharan, S.~Menon, A high-order adaptive {C}artesian cut-cell method
  for simulation of compressible viscous flow over immersed bodies, Journal of
  Computational Physics 321 (2016) 342--368.
\newblock \href {https://doi.org/10.1016/j.jcp.2016.05.050}
  {\path{doi:10.1016/j.jcp.2016.05.050}}.

\bibitem{kirkpatrick-cellLinking}
M.~Kirkpatrick, S.~Armfield, J.~Kent, A representation of curved boundaries for
  the solution of the {Navier–Stokes} equations on a staggered
  three-dimensional {C}artesian grid, Journal of Computational Physics 184~(1)
  (2003) 1--36.
\newblock \href {https://doi.org/10.1016/S0021-9991(02)00013-X}
  {\path{doi:10.1016/S0021-9991(02)00013-X}}.

\bibitem{cecere-cellLinking}
D.~Cecere, E.~Giacomazzi, An immersed volume method for large eddy simulation
  of compressible flows using a staggered-grid approach, Computer Methods in
  Applied Mechanics and Engineering 280 (2014) 1--27.
\newblock \href {https://doi.org/10.1016/j.cma.2014.07.018}
  {\path{doi:10.1016/j.cma.2014.07.018}}.

\bibitem{zeigelwanger-Pacman}
H.~Ziegelwanger, P.~Reiter, {The PAC-MAN model}: Benchmark case for linear
  acoustics in computational physics, Journal of Computational Physics 346
  (2017) 152--171.
\newblock \href {https://doi.org/10.1016/j.jcp.2017.06.018}
  {\path{doi:10.1016/j.jcp.2017.06.018}}.

\bibitem{sommariva-feketeQR}
A.~Sommariva, M.~Vianello, Computing approximate {F}ekete points by {QR}
  factorizations of {V}andermonde matrices, Computers \& Mathematics with
  Applications 57~(8) (2009) 1324--1336.
\newblock \href {https://doi.org/10.1016/j.camwa.2008.11.011}
  {\path{doi:10.1016/j.camwa.2008.11.011}}.

\bibitem{davis-momentFitting}
P.~J. Davis, A construction of nonnegative approximate quadratures, Mathematics
  of Computation 21~(100) (1967) 578--582.
\newblock \href {https://doi.org/10.2307/2005001} {\path{doi:10.2307/2005001}}.

\bibitem{hesthaven-TheDGBook}
J.~S. Hesthaven, T.~Warburton, Nodal Discontinuous {G}alerkin Methods:
  Algorithms, Analysis, and Applications.

\bibitem{chan-ESCurvilinear}
J.~Chan, L.~C. Wilcox, On discretely entropy stable weight-adjusted
  discontinuous {G}alerkin methods: curvilinear meshes, Journal of
  Computational Physics 378 (2019) 366--393.
\newblock \href {https://doi.org/10.1016/j.jcp.2018.11.010}
  {\path{doi:10.1016/j.jcp.2018.11.010}}.

\bibitem{chan-heterogenMedia}
J.~Chan, R.~J. Hewett, T.~Warburton, Weight-adjusted discontinuous {G}alerkin
  methods: Wave propagation in heterogeneous media, SIAM Journal on Scientific
  Computing 39~(6) (2017) A2935--A2961.
\newblock \href {https://doi.org/10.1137/16M1089186}
  {\path{doi:10.1137/16M1089186}}.

\bibitem{giuliani-weightedSRD}
A.~Giuliani, A.~Almgren, J.~Bell, M.~Berger, M.~{Henry de Frahan},
  D.~Rangarajan, A weighted state redistribution algorithm for embedded
  boundary grids, Journal of Computational Physics 464 (2022) 111305.
\newblock \href {https://doi.org/10.1016/j.jcp.2022.111305}
  {\path{doi:10.1016/j.jcp.2022.111305}}.

\bibitem{nordstrom-StableFiltering}
J.~Nordstr{\"o}m, A.~R. Winters, Stable filtering procedures for nodal
  discontinuous {G}alerkin methods, Journal of Scientific Computing 87~(1)
  (2021) 17.
\newblock \href {https://doi.org/10.1007/s10915-021-01434-x}
  {\path{doi:10.1007/s10915-021-01434-x}}.

\bibitem{lunquist-StableFiltering}
T.~Lundquist, J.~Nordstr{\"o}m, Stable and accurate filtering procedures,
  Journal of Scientific Computing 82~(1) (2020) 16.
\newblock \href {https://doi.org/10.1007/s10915-019-01116-9}
  {\path{doi:10.1007/s10915-019-01116-9}}.

\bibitem{jensens}
J.~L. W.~V. Jensen, Sur les fonctions convexes et les inégalités entre les
  valeurs moyennes, Acta Mathematica 30~(none) (1906) 175 -- 193.
\newblock \href {https://doi.org/10.1007/BF02418571}
  {\path{doi:10.1007/BF02418571}}.

\bibitem{tao-tunneling}
M.~Tao, C.~Batty, E.~Fiume, D.~I.~W. Levin, Mandoline: robust cut-cell
  generation for arbitrary triangle meshes, ACM Transactions on Graphics 38~(6)
  (nov 2019).
\newblock \href {https://doi.org/10.1145/3355089.3356543}
  {\path{doi:10.1145/3355089.3356543}}.

\bibitem{tsitouras-Tsit5}
C.~Tsitouras, {Runge–Kutta} pairs of order 5(4) satisfying only the first
  column simplifying assumption, Computers \& Mathematics with Applications
  62~(2) (2011) 770--775.
\newblock \href {https://doi.org/10.1016/j.camwa.2011.06.002}
  {\path{doi:10.1016/j.camwa.2011.06.002}}.

\bibitem{Julia-OrdinaryDiff}
C.~Rackauckas, Q.~Nie, {D}ifferentialequations.jl--a performant and
  feature-rich ecosystem for solving differential equations in {J}ulia, Journal
  of Open Research Software 5~(1) (2017) 15--15.
\newblock \href {https://doi.org/10.5334/jors.151}
  {\path{doi:10.5334/jors.151}}.

\bibitem{Julia-PathIntersections}
\href{https://github.com/cgt3/PathIntersections.jl}{{PathIntersections.jl}}.
\newline\urlprefix\url{https://github.com/cgt3/PathIntersections.jl}

\bibitem{Julia-StartUpDG}
\href{https://github.com/jlchan/StartUpDG.jl}{{StartUpDG.jl}}.
\newline\urlprefix\url{https://github.com/jlchan/StartUpDG.jl}

\bibitem{code-repo}
\href{https://github.com/cgt3/paper-2024-srd-wave}{Reproducability repository
  for simulation codes}.
\newline\urlprefix\url{https://github.com/cgt3/paper-2024-srd-wave}

\bibitem{Julia-ForwardDiff}
J.~Revels, M.~Lubin, T.~Papamarkou, Forward-mode automatic differentiation in
  {J}ulia (2016).
\newblock \href {http://arxiv.org/abs/1607.07892} {\path{arXiv:1607.07892}}.

\bibitem{saye-Quadrature}
R.~I. Saye, High-order quadrature methods for implicitly defined surfaces and
  volumes in hyperrectangles, SIAM Journal on Scientific Computing 37~(2)
  (2015) A993--A1019.
\newblock \href {https://doi.org/10.1137/140966290}
  {\path{doi:10.1137/140966290}}.

\bibitem{saye-QuadraturePolynomials}
R.~I. Saye, High-order quadrature on multi-component domains implicitly defined
  by multivariate polynomials, Journal of Computational Physics 448 (2022)
  110720.
\newblock \href {https://doi.org/10.1016/j.jcp.2021.110720}
  {\path{doi:10.1016/j.jcp.2021.110720}}.

\end{thebibliography}

\end{document}